\documentclass[a4paper,11pt]{amsart}

\usepackage{lineno,hyperref}

\usepackage{amsfonts,amsmath, amsthm, amscd, amssymb, graphicx, color, mathrsfs, stmaryrd}
\usepackage[utf8]{inputenc}

\def\findemo{\hfill \rule{6pt}{6pt}}
\newcommand{\conv}{\operatorname{conv}}
\newcommand{\dom}{\operatorname{dom}}
\newcommand{\diam}{\operatorname{diam}}
\newcommand{\Dz}{\operatorname{Dz}}
\newcommand{\epi}{\operatorname{epi}}

\newcommand{\dd}{\operatorname{d}}

\begin{document}

\newtheorem{thm}{Theorem}[section]
\newtheorem{theo}[thm]{Theorem}
\newtheorem{prop}[thm]{Proposition}
\newtheorem{coro}[thm]{Corollary}
\newtheorem{lema}[thm]{Lemma}
\newtheorem{defi}[thm]{Definition}
\newtheorem{ejem}[thm]{Example}
\newtheorem{rema}[thm]{Remark}
\newtheorem{fact}[thm]{Fact}
\newtheorem{open}[thm]{PROBLEM}

\title{On uniformly convex functions}

\author{G. Grelier, M. Raja}
\date{February 4, 2021}

\address{Departamento de Matem\'aticas, Universidad de Murcia, Campus de Espinardo, 30100 Espinardo, Murcia, Spain}

\thanks{The authors were supported by the Grants of Ministerio de Econom\'ia, Industria y Competitividad MTM2017-83262-C2-2-P; and Fundaci\'on S\'eneca Regi\'on de Murcia 20906/PI/18.}


\maketitle

\begin{abstract}
Non-convex functions that yet satisfy a condition of uniform convexity for non-close points can arise in discrete constructions. We prove that this sort of discrete uniform convexity is inherited by the convex envelope, which is the key 
to obtain other remarkable properties such as the coercivity. Our techniques allow to retrieve Enflo's uniformly convex renorming of super-reflexive Banach spaces  as the regularization of a raw function built from trees.
Among other applications, we provide a sharp estimation of the distance of a given function to the set of differences of Lipschitz convex functions.
Finally, we prove the equivalence of several possible ways to quantify the super weakly noncompactness of a convex subset of a Banach space.
\end{abstract}

\section{Introduction}

Along the paper, $(X, \| \cdot \|)$ will be a real Banach space and we will follow the standard notation that one can find in books such as \cite{BL, DGZ, banach, JL, LT}. However, dealing with real functions defined on $X$, if there is not specific hypothesis on the domain, we will follow the convention typical from Convex Analysis \cite{BV,Za2} that a function $f$ is defined everywhere and takes values in $\overline{{\Bbb R}}={\Bbb R} \cup\{-\infty,+\infty\}$. 
A function $f$ is said to be proper if $f(x) > -\infty$ for all $x \in X$, and $\dom(f):=\{x \in X: f(x) <+ \infty\} \not = \emptyset$.
In the following, all the functions are supposed to be proper. However, some operations performed on proper functions could lead to non-proper functions.
The class of lower semicontinuous convex proper functions on $X$ will be denoted $\Gamma(X)$. Note that, if nothing is
said on the contrary, all functions are supposed to be defined everywhere  in $X$ and proper.

\begin{defi}
Let $\varepsilon>0$.
A function $f: X \rightarrow \overline{{\Bbb R}}$ is said to be $\varepsilon$-uniformly convex if 
there is $\delta>0$ such that whenever $\|x-y\| \geq \varepsilon$, then  
$$f \left( \frac{x+y}{2} \right) \leq  \frac{f(x)+f(y)}{2}  - \delta. $$
The function is said to be uniformly convex if it is $\varepsilon$-uniformly convex for all $\varepsilon>0$.
\end{defi}

The suggestive name {\it discrete uniformly convex functions} applied to functions which are 
$\varepsilon$-uniformly convex for some $\varepsilon>0$ could be misleading here.
Clearly, a uniformly convex function is midpoint-convex, that is, the inequality $ f( \frac{x+y}{2} ) \leq  \frac{f(x)+f(y)}{2}$ holds whenever $x,y \in X$.
Therefore, a uniformly convex function is convex provided some regularity holds (e.g., if $f$ is lower semicontinuous).
The notion of uniform convexity for functions was introduced by Levitin and Polyak \cite{LePo}, and based on Clarkson's uniform convexity for normed spaces \cite{Cla}.
Since then, the properties of uniformly convex functions have been studied in several papers, notably  \cite{vlad1,vlad2,Za1,AP, BGHV, BV2}, the section 3.5 in Zalinescu's book \cite{Za2}, and part of chapter 5 in Borwein-Vaderwerff's book \cite{BV} devoted to them. In relation to the standard theory, let us point out the notion of modulus of uniform convexity 
$$ \delta_f(\varepsilon) = \inf \left\{  \frac{f(x)+f(y)}{2} - f \left(\frac{x+y}{2} \right): x,y \in \dom(f), \|x-y\| \geq \varepsilon \right\} .$$
Note that $ \delta_f$ could take negative values unless $f$ is supposed to be (midpoint-) convex.
Analogously, it is possible to define $\varepsilon$-uniformly concave functions, however it will not be necessary to treat them here because all the theory extends trivially.\\

In this paper, we are focused in $\varepsilon$-uniformly convex functions for a fixed $\varepsilon>0$ for which  
{\bf the usual convexity assumption is not longer assumed}. 
That is the main issue we have to deal with here and the reason to do it 
is that non-convex $\varepsilon$-uniformly convex functions
may arise in relation with some discrete constructions, starting from trees or barely convex sets. 
Nonetheless, $\varepsilon$-uniformly convex functions have nice properties.
Along the paper, $\breve{f}$ will denote the lower semicontinous convex envelope of a function $f$ (also denoted 
$\overline{\conv}(f)$ in some references). The next result shows the global behaviour of $\varepsilon$-uniformly convex functions and the relative stability of minimizers by linear perturbations.

\begin{theo}\label{main-coercive}
Let $f$ be an $\varepsilon$-uniformly convex function such that $\breve{f}$ is proper. Then $f$ is bounded below and coercive, more precisely we have   
$$ \liminf_{\|x\| \rightarrow +\infty} \frac{f(x)}{~ \|x\|^2} > 0 .$$
Moreover, for any $\varepsilon'>\varepsilon$ there exist $\delta, \eta>0$ such that if given $x_0^* \in X^*$ and $x_0 \in X$ with
$$f(x_0)+x_0^*(x_0) < \inf (f+x_0^*)+\delta,$$
and $x^* \in X^*$ such that $\|x^* -x_0^*\| < \eta$ and $x \in X$ that minimizes $f+x^*$, then $\|x-x_0\| \leq \varepsilon'$.
The existence of such minimizer pair $(x,x^*)$ is guaranteed if $f=\breve{f}$.
\end{theo}

The proof of the former result relies in the possibility of ``making convex'' an $\varepsilon$-uniformly convex function without loosing the $\varepsilon$-uniformly convexity. We will say that a function $f$ is $\varepsilon^+$-uniformly convex if it is $\varepsilon'$-uniformly convex for every $\varepsilon'>\varepsilon$. We have the following result.

\begin{theo}\label{main-convexify}
Let $f$ be $\varepsilon$-uniformly convex and assume that $\breve{f}$ is proper. Then
$\breve{f}$ is $\varepsilon^+$-uniformly convex.
\end{theo}

Simple examples, such as Example \ref{ejemb}, show that the $\varepsilon$-uniformly convexity of 
$f$ does not guarantee that $\breve{f}$ would be proper. In order to fulfil that requirement in terms of $f$, we
direct the reader to Corollary \ref{coro-proper}.
Suppose now that we already have a proper lower semicontinuous convex and $\varepsilon$-uniformly convex function $f$. We wonder if we could ``upgrade'' $f$ to a new function sharing those properties and, besides, being locally Lipschitz (global Lipschitzness is not allowed for uniformly convex functions). In that sense, we have the following result.

\begin{theo}\label{main-renorming}
Let $f \in \Gamma(X)$ be $\varepsilon$-uniformly convex. Then there exists an equivalent norm $|\!|\!|  \cdot |\!|\!| $ on $X$ such that the function $x \rightarrow |\!|\!|  x |\!|\!| ^2 $ is $\varepsilon^+$-uniformly convex on the subsets of $\dom(f)$ where $f$ is bounded above.
Moreover, the norm  $|\!|\!|  \cdot |\!|\!| $ can be taken as close to $\| \cdot \|$ as we wish.
\end{theo}

We want to point out that in the previous theorem we get $\varepsilon^+$ even in the case that the function $f$ in the hypothesis should be convex. If $f$ were uniformly convex, then a series of $\varepsilon$-uniformly convex norms for different $\varepsilon$'s going to $0$ would produce an equivalent norm whose square is uniformly convex on bounded subsets of $\dom(f)$.\\

It turns out that supporting a convex continuous $\varepsilon$-uniformly convex function is actually a geometrical-topological property of the domain.
It is known that a Banach space admits a uniformly convex function bounded on bounded sets if and only if it is super-reflexive. The second named author proved in \cite{raja2} that a closed convex bounded set admits a bounded continuous uniformly convex function if and only if it is {\it super weakly compact} (SWC for short). We will give the actual definition of SWC set in Section 6, however we can provide an alternative one on a provisional basis: a bounded closed convex set is SCW if and only if  for all $\varepsilon>0$ there is $N_\varepsilon \in {\Bbb N}$ such that the height of any $\varepsilon$-separated dyadic tree is bound by $N_\varepsilon$. Recall that 
a {\it dyadic tree} of height $n \in {\Bbb N}$ is a set of the form $\{ x_s: |s| \leq n \}$,
indexed by finite sequences $s \in \bigcup_{k=0}^n \{0,1\}^k$ of length $|s| \leq n$, such that
$x_s=2^{-1} ( x_{s\frown 0} +  x_{s\frown 1})$ for every $|s| < n$, where
$\{0,1\}^0:=\{\emptyset\}$ indexes the root $x_\emptyset$
and the symbol ``$\frown$'' stands for concatenation.
We say that a dyadic tree $\{ x_s: |s| \leq n \}$ is
{\it $\varepsilon$-separated} if $\| x_{s\frown 0} - x_{s\frown 1} \| \geq \varepsilon $ for every $|s|<n$.\\

Our techniques allow us to give a very precise quantitative version of the relation between containment of separated trees and supporting a uniformly convex function for a set.

\begin{theo}\label{main-swc}
Let $C \subset X$ be a closed bounded convex set. Then these two numbers coincide:
\begin{itemize}
\item[$(a)$] the infimum of the $\varepsilon>0$ such that there is a common bound for the heights of all the $\varepsilon$-separated dyadic trees;
\item[$(b)$]  the infimum of the $\varepsilon>0$ such that there is a bounded $\varepsilon$-uniformly convex (and convex, Lipschitz\dots) function defined on $C$.
\end{itemize}
\end{theo}

As we will see later, the quantities given by the previous theorem can be used as measures of super weak noncompactness. 
Note that the combination of Theorem \ref{main-swc} and Theorem \ref{renorming} applied to the ball of a Banach space produces yields the famous Enflo's renorming theorem of super-reflexive spaces.
At this point, we want to stress that we barely get Enflo's but not Pisier's, see \cite{DGZ} for instance, because we are mainly focused on ``$\varepsilon$'' (the separation of dyadic trees) instead of ``$\delta$'' (the quality of the modulus of convexity).
The importance of Enflo's result motivated us to offer to the reader a more direct proof based on our arguments.\\ 

A couple of comments on the contents of this paper. We will consider the more general notion of  {\it $\varepsilon$-uniformly convexity with respect to a metric $d$}, instead of the norm, or even a {\it pseudometric}. Namely, let $d$ be a pseudometric 
defined on $\dom(f)$ 
(that we will always assume uniformly continuous with respect to $\| \cdot \|$ by technical reasons). Given $\varepsilon>0$, we say that $f$ is
$\varepsilon$-uniformly convex with respect to $d$ if 
there is $\delta>0$ such that if $d(x,y) \geq \varepsilon$ then  
$$f \left(\frac{x+y}{2} \right) \leq  \frac{f(x)+f(y)}{2}  - \delta $$
(the modulus $\delta_f$ is defined likewise). With this definition Theorem~\ref{main-convexify}, 
Theorem~\ref{main-renorming} and Theorem~\ref{main-swc} are still true provided that $\dom(f)$ is bounded.
It is known that the dual notion of uniform convexity is the {\it uniform smoothness} \cite{AP,Za2, BGHV}, however, 
we will not discuss Fenchel duality here for $\varepsilon$-uniformly convex  functions. That will be eventually done in a subsequent paper.\\

The structure of the paper is the following. The second section deals with basic properties of $\varepsilon$-uniformly convex  and $\varepsilon$-uniformly quasi-convex functions, mostly under the hypothesis of convexity. A few examples are given to show that the definitions do not guarantee some additional nice properties. The third section is devoted to the proof of Theorem \ref{main-convexify} that will allow the reduction to the convex case of other results. The construction of uniformly convex functions form scratch (trees and sets) is done in the fourth section. The fifth section treats general properties of $\varepsilon $-uniformly convex functions and the possibility of adding more properties like Lipschitzness or homogeneity (renorming). We also prove an estimation of the approximation by differences of convex functions. In the sixth section, we prove the equivalence of several measures of super weak noncompactness for convex sets. We also propose a  measure of super weak noncompactness for bounded sets and we study its behaviour by convex hulls. 
In the last section we will sketch an understandable proof of Enflo's uniformly convex renorming of super-reflexive spaces theorem based on the ideas exposed along the paper.

\section{Basic properties and examples}

We will discuss in this section results of almost arithmetical nature.
The first proposition contains some easy facts whose proof is left to the reader.

\begin{prop}\label{basics}
Let $\varepsilon>0$ and let $f$ be an $\varepsilon$-uniformly convex function. Then:
\begin{enumerate}
    \item If $g$ is convex, then $f+g$ is $\varepsilon$-uniformly convex with $\delta_{f+g} \geq \delta_f$.
    \item The supremum of finitely many $\varepsilon$-convex functions is $\varepsilon$-convex too.
    \item If $f \geq 0$, then $f^2$ is $\varepsilon$-uniformly convex.
    \item The lower semicontinuous envelope of $f$ is $\varepsilon$-uniformly convex.
\end{enumerate}
\end{prop}

Recall that the infimal convolution of two functions $f,g$ is defined as 
$$ (f \, \square \, g)(x)= \inf \{f(x-y)+g(y): y \in X \}, ~\mbox{for}~x \in X. $$

\begin{prop}
Let $f_1,f_2$ be two convex functions such that $f_1$ is $\varepsilon_1$-uniformly convex and 
$f_2$ is $\varepsilon_2$-uniformly convex for $\varepsilon_1,\varepsilon_2>0$. Then  $f_1 \, \square \, f_2$ is $(\varepsilon_1+\varepsilon_2)$-uniformly convex with modulus $\min\{ \delta_{f_1}(\varepsilon_1),  \delta_{f_2}(\varepsilon_2) \}$.
\end{prop}

\noindent
{\bf Proof.} Given $x_1,x_2 \in \dom(f_1 \, \square \, f_2)=  \dom(f_1)+ \dom(f_2)$ with $\|x_1-x_2\| \geq \varepsilon_1+\varepsilon_2$ and $\eta>0$ we may find $y_1,y_2 \in \dom(f_2)$ such that 
$$  f_1(x_1-y_1) + f_2(y_1) < (f_1 \, \square \, f_2)(x_1) + \eta ,$$
$$  f_1(x_2-y_2) + f_2(y_2) < (f_1 \, \square \, f_2)(x_2) + \eta .$$
We have
$$ \| (x_1-y_1) - (x_2-y_2) \| + \| y_1 -y_2 \| \geq \|x_1-x_2\| \geq \varepsilon_1+\varepsilon_2. $$
Therefore, one of the inequalities either
$$ \| (x_1-y_1) - (x_2-y_2) \| \geq \varepsilon_1 ~~\mbox{or}~~ \|y_1-y_2\| \geq \varepsilon_2 $$
holds. Assume the first one does (the other case is similar)
$$  (f_1 \, \square \, f_2) \left( \frac{x_1+x_2}{2} \right) \leq f_1 \left(\frac{x_1+x_2}{2} - \frac{y_1+y_2}{2} \right) +
f_2 \left(  \frac{y_1+y_2}{2} \right)$$
$$ \leq  \frac{f_1(x_1-y_1) + f_1(x_2-y_2) }{2} - \delta_{f_1} (\varepsilon_1) + \frac{f_2(y_1)+f_2(y_2)}{2}$$
$$ \leq \frac{ (f_1 \, \square \, f_2)(x_1) + (f_1 \, \square \, f_2)(x_2) }{2} - \delta_{f_1}(\varepsilon_1) + \eta $$
which implies the statement as $\eta>0$ was arbitrary.\findemo\\

Now we will discuss some properties of the modulus of uniform convexity in the classic case, that is, when the function is assumed to be also convex.
The following property can be deduced easily with the help of a picture.

\begin{prop}\label{prevprop}
Let $f$ be convex and $\varepsilon>0$. Then 
$$ (1-t) f(x) + t f(y) - f((1-t)x+ty) \geq 2\delta_f(\varepsilon) \min\{t,1-t\}  $$
whenever $x,y \in \dom(f)$, $\|x-y\| \geq \varepsilon$ and $t \in [0,1]$.
\end{prop}

\noindent
{\bf Proof.} Without loss of generality we may assume $t \in [0,1/2]$ so $t = \min\{t,1-t\}$. Note now that
$$ (1-t)x+ty = (1-2t)x+2t\frac{x+y}{2} .$$
By convexity of $f$ we have
$$ f((1-t)x+ty) \leq (1-2t)f(x)+2tf \left( \frac{x+y}{2} \right) $$
$$\leq 
(1-2t)f(x)+2t \left( \frac{f(x)+f(y)}{2} -\delta_f(\varepsilon) \right)
= (1-t)f(x)+tf(y) - 2t\delta_f(\varepsilon) $$
as wished.\findemo\\

\begin{figure}
\centerline{\includegraphics[height=40mm]{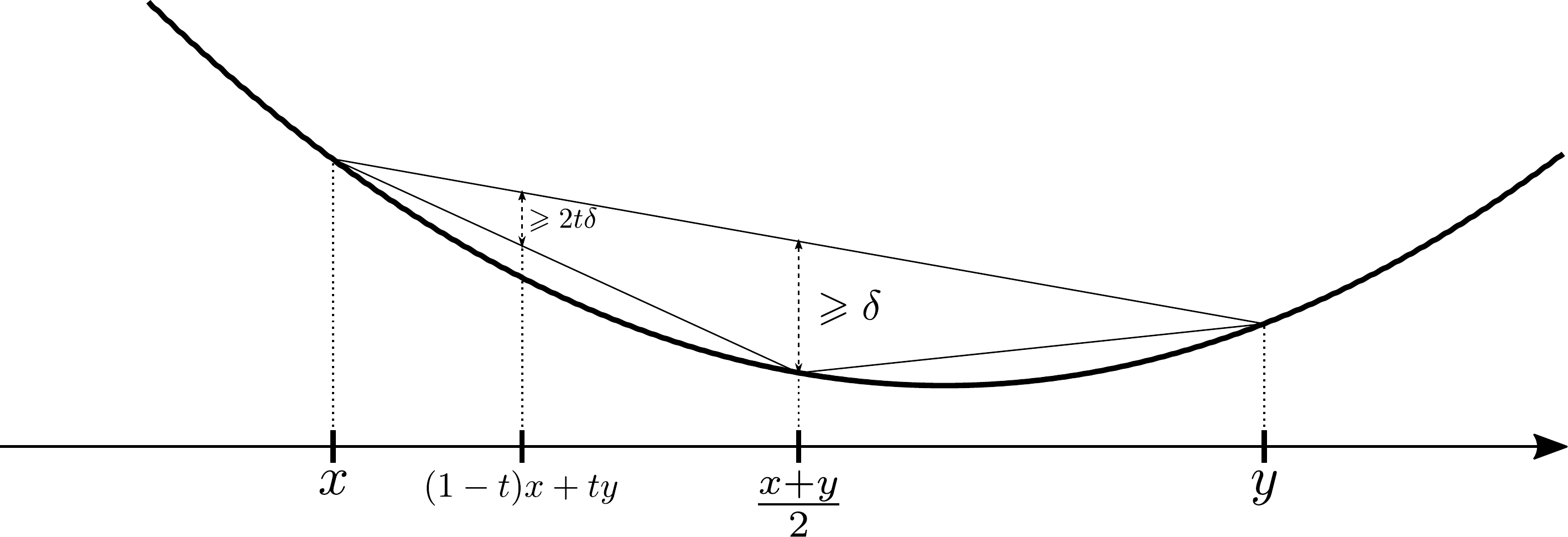}}
\caption{Geometric interpretation of Proposition \ref{prevprop}.}
\end{figure}

The {\it gage of uniform convexity} is introduced in \cite{vlad1} (see also \cite[p. 203]{Za2}) for convex function as 
$$ p_f(\varepsilon) = \inf\left\{ \frac{(1-t)f(x)+tf(y)-f((1-t)x+ty)}{t(1-t)}: 0 <t<1, \|x-y\| \geq \varepsilon\right\} .$$

\begin{coro}
For any convex function $f$ defined on $X$ and $\varepsilon>0$, we have
$$ 2 \delta_f(\varepsilon) \leq p_f(\varepsilon) \leq 4 \delta_f(\varepsilon).  $$
\end{coro}

\noindent
{\bf Proof.} 
The first inequality is a consequence of Proposition \ref{prevprop} together with the fact that $\min\{t,1-t\} \geq t(1-t)$, for $t \in {\Bbb R}$. The second inequality follows just taking $t=1/2$.\findemo\\

Therefore, for convex functions, $\varepsilon$-uniformly convexity can be expressed as $p_f(\varepsilon)>0$. The gage of uniform convexity has the following remarkable property 
$$ p_f(\lambda \varepsilon) \geq \lambda^2 p_f(\varepsilon) $$
whenever $\varepsilon \geq 0$ and $\lambda \geq 1$, see \cite[Proposition 3.5.1]{Za2} and note that the proof does not requiere the uniform convexity of $f$. In particular $\varepsilon \rightarrow \varepsilon^{-2} p_f (\varepsilon)$ is a non decreasing function.\\

Now we will discuss some examples showing the limitations of the notions we are dealing with.

\begin{figure}\label{fig2}
\centerline{\includegraphics[width=10cm]{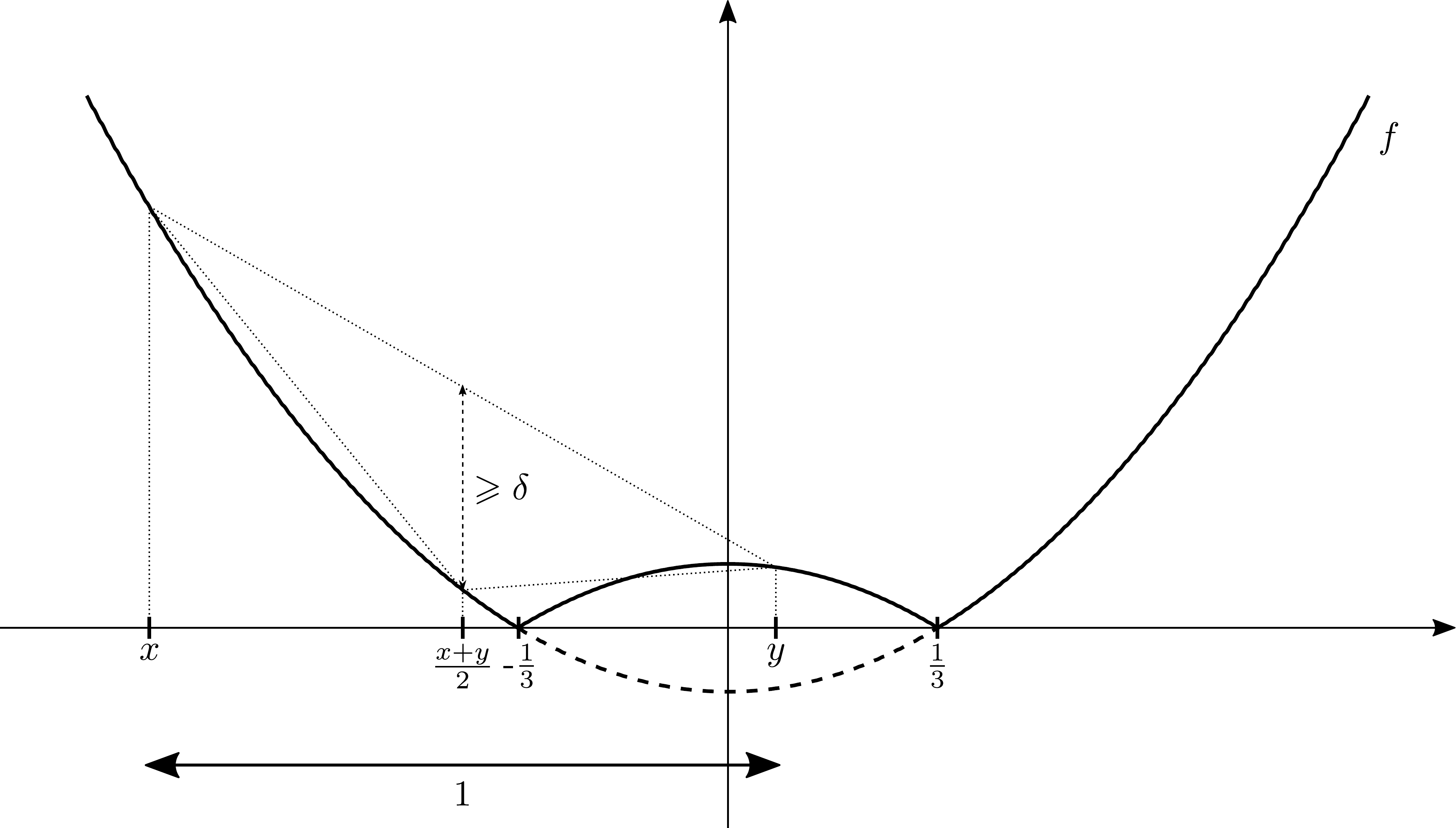}}
\caption{Example \ref{ejema}.}
\end{figure}

\begin{ejem}\label{ejema}
$f(x)=|x^2-1/9|$ is a continuous nonconvex $1$-uniformly convex function on ${\Bbb R}$.
\end{ejem}

\noindent
{\bf Proof.} This can be deduced by inspection of the drawing (see Figure \ref{fig2}). A more detailed computation shows that $\delta=1/36$.\findemo

\begin{ejem}\label{ejemb}
A (proper) $\varepsilon$-uniformly convex function may have a non-proper lower semicontinuous convex envelope.
\end{ejem}

\noindent
{\bf Proof.} Take a function $f$ which is finite and  unbounded below on $B(0,\varepsilon/3)$ and takes the value $+\infty$ outside. By the very definition, $f$ is
$\varepsilon$-uniformly convex and necessarily $\breve{f}=-\infty$ on $B(0,\varepsilon/3)$.\findemo

\begin{ejem}
A uniformly convex continuous function taking finite values which is unbounded on a bounded convex closed set.
\end{ejem}

\noindent
{\bf Proof.} The function will be defined on $\ell^2$. Firstly note that $\|x\|^2$ is uniformly convex. 
Consider the convex function
$h:{\Bbb R} \rightarrow {\Bbb R}$ defined by
$$ h(t) = \max\{0, t -1/2, -t -1/2\} .$$
The series $g(x)=\sum_{n=1}^\infty n h(x_n)$, for $x=(x_n) \in \ell^2$, defines a convex continuous function. Indeed, at each point, only a finite number of summands can be positive at once. The continuity comes from the fact that the same is true on 
any ball of radius less than $1$. Now, the function $f(x)=\|x\|^2+g(x)$ is continuous, unbounded on $B_{\ell^2}$ and, by Proposition \ref{basics}, it is also uniformly convex.\findemo\\

The following notions will be useful in relation with $\varepsilon$-uniform convexity. 

\begin{defi}
Let $f:X \rightarrow \overline{{\Bbb R}}$ be a function. Then $f$ is said to be:
\begin{enumerate}
\item quasi-convex if 
$$ f(\lambda x + (1-\lambda) y) \leq \max\{f(x),f(y)\} $$
for every $x,y \in X$ and $\lambda \in [0,1]$.
\item $\varepsilon$-uniformly quasi-convex if, for a given $\varepsilon>0$, there is some $\delta>0$ such that
$$ f \left( \frac{x+y}{2} \right) \leq \max\{f(x),f(y)\} - \delta $$
whenever $x,y \in X$ with $\|x-y\| \geq \varepsilon$ (or $d(x,y) \geq \varepsilon$ for a pseudometric $d$).
\item uniformly quasi-convex if it is $\varepsilon$-uniformly quasi-convex for every $\varepsilon>0$.
\end{enumerate}
\end{defi}

Whereas the notion of quasi-convexity is well known, our definition of uniform quasi-convexity is weaker than the one given in \cite{vlad2}. As with convexity, the midpoint version does not implies the ``$\lambda$-version'' unless some regularity (e.g. lower semicontinuity) is assumed. The following result shows one relation between the quantified versions of uniform convexity and uniform quasi-convexity for functions. 

\begin{prop}\label{square-quasi}
Let $\varepsilon>0$ and let $f \geq 0$ be a convex and $\varepsilon$-uniformly quasi-convex function. 
Then $f^2$ is $\varepsilon$-uniformly convex.
\end{prop}

\noindent
{\bf Proof.}
The following inequality can be checked easily: if for some real numbers $a,b,c$ we have $a+b \geq 2c \geq 0$ then
\begin{equation}
\label{inequality}
\left( \frac{a+b}{2} - c \right)^2+
\left( \frac{a-b}{2} \right)^2 \leq 
\frac{a^2+b^2}{2} - c^2.
\end{equation}
Assume $\|x-y\| \geq \varepsilon$ and let $\delta>0$ be given by the definition of $\varepsilon$-uniform quasi-convexity. 
If $|f(x)-f(y)| > \delta$ the previous inequality implies 
$$ \frac{f(x)^2+f(y)^2}{2} - f \left( \frac{x+y}{2} \right)^2 \geq \frac{\delta^2}{4} .$$
On the other hand, if $|f(x)-f(y)| \leq \delta$ then
$$ f \left( \frac{x+y}{2} \right) \leq \max\{f(x),f(y)\} -\delta \leq \frac{f(x)+f(y)}{2} - \frac{\delta}{2} $$
and thus
$$ \frac{f(x)^2+f(y)^2}{2} - f \left( \frac{x+y}{2} \right)^2 \geq \frac{\delta^2}{4} $$
using again the inequality (\ref{inequality}).\findemo 

\begin{ejem}
A uniformly quasi-convex non-convex (concave) function. 
\end{ejem}

\noindent
{\bf Proof.} Take $f(x)=x$ for $x<0$, and $f(x)=x/2$ for $x \geq 0$.\findemo

\section{Convexifying the $\varepsilon$-uniform convexity}

In order to cover previous developments around finite dentability \cite{raja} we will consider uniformly convex functions with respect to a pseudometric $d$ defined on the domain of $f$. The norm of the Banach space will still play an important role and we requiere that $d$ be uniformly continuous with respect to the norm. Therefore, along this section we will assume that 
$\varepsilon$-uniform convexity refers to $d$. We will refer as $d$-diameter of a subset in $X \times {\Bbb R}$ the diameter with respect to $d$ of the projection of the set onto $X$. Let $\varpi$ be the modulus of uniform continuity (the standard symbol is ``$\omega$'' but we are using it as the first countable ordinal later), that is, the following inequality holds$$d(x,y) \leq \varpi(\|x-y\|)$$
and $\lim_{t \rightarrow 0^+} \varpi(t)=0.$

\begin{prop}\label{criterio-slice}
Let $f$ be a function and let $\varepsilon>0$. Then
\begin{enumerate}
\item If $f$ is $\varepsilon$-uniformly convex then every slice of $\epi(f)$ disjoint from 
$\epi(f+\delta_{f}(\varepsilon))$ has $d$-diameter less than $\varepsilon$.
\item If $f \in \Gamma(X)$  and there is $\delta>0$ such that every slice of $\epi(f)$ disjoint from 
$\epi(f+\delta)$ has $d$-diameter less than $\varepsilon$ then $f$ is $\varepsilon$-uniformly convex with modulus 
$\delta_{f}(\varepsilon) \geq \delta/2$.
\end{enumerate}
\end{prop}

\noindent
{\bf Proof.} For the first statement, assume that $(x,f(x)),(y,f(y))$ belong to such a slice. The separation from 
$\epi(f+ \delta_{f}(\varepsilon) )$ implies
$$ \frac{f(x)+f(y)}{2} < f \left( \frac{x+y}{2} \right) + \delta_{f}(\varepsilon)  $$
and so $d(x,y)<\varepsilon$. 
On the other hand, let $\delta>0$ as in statement 2 and take $x,y \in X$ such that the following inequality holds
$$ \frac{f(x)+f(y)}{2} - f \left( \frac{x+y}{2} \right) <\frac{\delta}{2}.$$
It implies that $\left( \frac{x+y}{2}, \frac{f(x)+f(y)}{2} \right)$ does not belong to $\epi(f+\delta/2)$. 
We may take an affine function $h$ such that $h < f + \delta/2$ and $h(\frac{x+y}{2}) > \frac{f(x)+f(y)}{2}$. It is evident that either $f(x) < h(x)$ or $f(y) < h(y)$. We may assume without loss of generality that the first inequality holds as the scenario is symmetric for $x$ and $y$. 
Now we have 
$$ f(y) < 2h \left( \frac{x+y}{2} \right) - f(x) = h(x) + h(y) -f(x) < h(y) + \frac{\delta}{2}  .$$
That implies both $(x,f(x))$ and $(y,f(y))$ belong to the slice defined by $h+\delta/2$ 
$$ S= \{ (x,t) \in \epi(f): t < h(x)+\delta/2 \} .$$
By our choices, we have $S \cap \epi(f+\delta)=\emptyset$ and thus 
$d(x,y) < \varepsilon$ by the hypothesis. We deduce in this way that $\delta/2 \leq \delta_{f}(\varepsilon)$.\findemo

\begin{coro}\label{convex-delta}
Let $\varepsilon>0$ and let $f$ be a convex and $\varepsilon$-uniformly convex function. Then
$$ f(x) \leq \sum_{k=1}^n \lambda_k f(x_k)- \delta_{f}(\varepsilon) $$
whenever  $x,x_1,\dots,x_n \in \dom(f)$ satisfy that $d(x,x_k) \geq \varepsilon$ and
$x=\sum_{k=1}^n \lambda_k x_k$ with $\lambda_k \geq 0$ and $\sum_{k=1}^n \lambda_k =1$. 
\end{coro}

\noindent
{\bf Proof.} If the inequality does not hold, then $ \left(x,  \sum_{k=1}^n \lambda_k f(x_k) \right)$ does not belong to $\epi(f+\delta_{f}(\varepsilon))$
so it can be separated from that set with a slice. Necessarily, one of the points $(x_k,f(x_k))$ belongs to the slice. That implies that the $d$-diameter of the slice is at least $\varepsilon$ which contradicts the previous proposition.\findemo\\

The following result is based on the techniques of the geometrical study of the Radon--Nikodym property (RNP), see \cite{bourgin}.
Note that the technique works only on bounded domains.

\begin{lema}
Let $\varepsilon>0$ and let $f$ be a bounded below function with bounded domain. 
Let $m>0$ be an upper bound for the norm diameter of $dom(f)$ and let 
$\tau>0$ be such that $\tau/m<1$. Assume that the set 
$$\{x \in X: f(x) < \inf f +\delta\}$$ 
has $d$-diameter less than $\varepsilon$. Then the set 
$$\{x \in X: \breve{f}(x) < \inf f + \delta \tau/m\}$$ 
has $d$-diameter less than $\varepsilon+2\varpi(\tau)$.
\end{lema}

\noindent
{\bf Proof.} Consider the sets
$$A = \{(x,r) \in X\times {\Bbb R}: f(x) \leq r < \inf f +\delta\};~~$$
$$B = \{(x,r) \in X\times {\Bbb R}:  r \geq \inf f +\delta, f(x) \leq r\}. $$
Note that the epigraph of $f$ is $A \cup B$. Consider their closed convex hulls $\breve{A}=\overline{\conv}(A)$ and
$\breve{B}=\overline{\conv}(B)$ and note that $\conv(\breve{A} \cup \breve{B})$ is dense in the epigraph of 
$\breve{f}$. 
Assume that $(x,r) \in \conv(\breve{A} \cup \breve{B})$ and $r <  \inf f +\delta\tau/m$. There is $\lambda \in [0,1]$ such that $(x,r)=\lambda(y,t) + (1-\lambda)(z,s)$ where $(y,t) \in \breve{A}$ and
$(z,s) \in \breve{B}$. 
The condition $\lambda t + (1-\lambda) s <  \inf f +\delta\tau/m$ implies $1-\lambda< \tau/m$.
Indeed, suppose $1-\lambda \geq \tau/m$. As $s \geq  \inf f + \delta $, then 
$$(1-\lambda) s \geq  (1-\lambda) \inf f + \delta \tau/m $$ 
On the other hand, $ \lambda t \geq \lambda \inf f$. Adding these inequalities we get $\lambda t + (1-\lambda) s \geq  \inf f +\delta\tau/m$, which contradicts the assumption.
Therefore
$$ \|x-y\| = \| (\lambda -1) y + (1-\lambda) z\|= (1-\lambda) \|y-z\| < \tau . $$
In order to estimate the $d$-diameter of 
$$ S = \{ (x,t): x \in X,  \breve{f}(x) \leq t < \inf f + \delta\tau/m\}, $$
we may consider only points on the dense set $S \cap  \conv(\breve{A} \cup \breve{B})$.
Therefore, consider 
$(x_1,r_1), (x_2,r_2) \in \conv(\breve{A} \cup \breve{B})$ with $r_1,r_2 < \inf f +\delta\tau/m$.
The convex decomposition above shows that for some $\lambda_1,\lambda_2 \in [0,1]$ and points $(y_1,t_1),(y_2,t_2)  \in \breve{A}$ and $(z_1,s_1), (z_2,s_2) \in \breve{B}$
we have
$$(x_1,r_1)=\lambda_1(y_1,t_1) + (1-\lambda_1)(z_1,s_1),$$ 
$$(x_2,r_2)=\lambda_2(y_2,t_2) + (1-\lambda_2)(z_2,s_2).$$ 
By the previous estimations, we have $\|x_1-y_1\|, \|x_2-y_2\| \leq \tau$, which implies that  $d(x_1,y_1),d(x_2,y_2) \leq \varpi(\tau)$, and thus, by the assumption on $A$,  
$$ d(x_1,x_2) \leq d(x_1,y_1)+ d(y_1,y_2)+ d(x_2,y_2) \leq \varepsilon+ 2\varpi(\tau) $$
as desired.\findemo\\

We will deal now with the transfer of the $\varepsilon$-uniform convexity property to the lower semicontinuous convex envelope. Note that $\varepsilon$-uniform convexity is referred to a uniformly continuous pseudometric, however we require the hypothesis that the domain be norm bounded.

\begin{theo}\label{convexify}
Let $\varepsilon>0$ and let $f$ be a bounded below $\varepsilon$-uniformly convex function with bounded domain . 
Then $\breve{f}$ is  $\varepsilon^+$-uniformly convex and given $\varepsilon'>\varepsilon$, the modulus of convexity  
$\delta_{\breve{f}}(\varepsilon')$ depends only on $\varepsilon'$, $\delta_f(\varepsilon)$, $\varpi$ and the norm diameter of 
$\dom(f)$.
\end{theo}

\noindent
{\bf Proof.} 
Let $m$ an upper bound for the diameter of $\dom(f)$ and $\delta>0$ the parameter given by the definition of $\varepsilon$-uniform convexity. Take $\tau>0$ such that $\tau/m<1$. We will estimate the $d$-diameter of any slice of $\epi(\breve{f})$ not meeting $\epi(\breve{f}+\delta\tau /m)$. 
Suppose that the slice is given by $x^* \in X^*$. Note that the estimation of the $d$-diameter of the slice we 
need is equivalent to the same for an horizontal slice of 
$\epi(\breve{f}-x^*)$ not meeting  $\epi(\breve{f}-x^*+ \delta\tau/m)$, 
which is the same as taking the points of 
$\epi(\breve{f}-x^*)$ whose scalar coordinate is less than $\inf (\breve{f}-x^*)+ \delta\tau/m$.
Since $\breve{f}-x^*$ equals the convex envelope of the function $f-x^*$, which is $\varepsilon$-uniformly convex with parameter $\delta$, the set 
$$\{x \in X: f(x)-x^*(x) < \inf (f-x^*) +\delta\}$$ 
has diameter less than $\varepsilon$ by Proposition \ref{criterio-slice}. The previous lemma applies to get that
$$\{x \in X: \breve{f}(x)-x^*(x) < \inf ( f -x^* ) + \delta\tau/m\}$$ 
has diameter less than $\varepsilon + 2\varpi(\tau)$. 
Thanks to Proposition \ref{criterio-slice}, it follows that $ \breve{f}$ is 
$\varepsilon + 2\varpi(\tau)$-uniformly convex. Given $\varepsilon'>\varepsilon$, we only have to set 
 $\tau>0$ such that $2\varpi(\tau)<\varepsilon'-\varepsilon$.\findemo\\

The following result is the key to deal with unbounded domains.

\begin{prop}\label{bounded-reduction}
Let $\varepsilon>0$ and let $f$ be an $\varepsilon$-uniformly convex function such that $\breve{f}$ is proper. Then the value of $\breve{f}(x)$ for $x \in \dom(f)$ depends only on the set of values $\{f(y): \|y-x\| <\varepsilon \}$. Namely,  if  
$g$ is the function defined by $g(y)=f(y)$ if $\|y-x\| <\varepsilon$ and $g(x)=+\infty$ otherwise, 
then $\breve{f}(x)=\breve{g}(x)$.
\end{prop}

\noindent
{\bf Proof.}
Let us roughly explain the idea of the proof before going into details. A priori, the computation of $\breve{f}(x)$ may involve values of $f$ at points arbitrarily far away from $x$. Namely, $(x,\breve{f}(x))$ can be approximated by a convex combination of points of the form $(x_k,f(x_k))$. As we want the points $x_k$ to be close to $x$, we will describe an algorithm that will modify the set $\{x_k\}$ by the substitution of one (or several points) at each step until the resulting set is contained in $B(x,\varepsilon)$. The algorithm consist in switching a farthest point $x_i$ by the middle point between it and an ``opposite point'' $x_j$ which is not farther from $x$ as $x_i$ is. If $d(x_i,x_j) \geq \varepsilon$, the $\varepsilon$-uniform convexity of $f$ will imply that we do not loose information about $\breve{f}(x)$ when switching $x_i$ by $(x_i+x_j)/2$. Once, $x_i$ has disappear from the the set, we choose a new farthest point and start over. 
Actually, the method bring the points closer to $x$ with respect to a prefixed direction $x^* \in X^*$. The repetition of the algorithm with several directions will eventually finish with the modified set of points contained into $B(x,\varepsilon)$.
Now we will resume the proof.\\
The definition of $\breve{f}$ implies that the following set 
$$ \left\{ (x,t): t \geq  \sum_{k=1}^n \lambda_k f(x_k) ~~ \mbox{with}~~
x=\sum_{k=1}^n \lambda_k x_k ~~ \mbox{a convex combination}  \right\} $$
is dense in $\epi(\breve{f})$. Fix $x \in \dom(f)$ and suppose $x=\sum_{k=1}^n \lambda_k x_k$ is a convex combination.
Now, we are going to describe the announced algorithm that will transform the set of points $S=\{x_1,\dots,x_n\}$ into a set $S'=\{x'_1,\dots,x'_{n'}\}
\subset B(x,\varepsilon)$ such that still we have $\sum_{k=1}^{n'} \lambda'_k x'_k=x$, where
$\sum_{k=1}^{n'} \lambda'_k x'_k=x$ with $\lambda'_k \geq 0$, $\sum_{k=1}^{n'} \lambda'_k=1$, and
$$ \sum_{k=1}^{n'} \lambda'_k f(x'_k) \leq  \sum_{k=1}^n \lambda_k f(x_k) .$$
In order to do that, without loss of generality, we may assume $x=0$. Fix $x^* \in S_{X^*}$. Let $a= \sup x^*(S) \geq 0$ and $b=-\inf x^*(S)$. As $a$ and $b$ can be exchanged just taking $-x^*$ instead, without loss of generality we may 
assume $a\geq \max\{b,\varepsilon\}$.
Also, without loss of generality, we may assume $x^*(x_1)=a$. Since $x_1$ is the farther point (with respect to $x^*$), its ``mass'' $\lambda_1$ compensates with masses on the side $x^* \leq 0$. Suppose firstly that $x^*(x_2) \leq 0$ and $\lambda_2 \geq \lambda_1$. We have $\|x_1-x_2\| \geq \varepsilon$.
We claim that it is possible to switch $x_1$ by $x_1'=(x_1+x_2)/2$. Indeed, 
$$ 2\lambda_1 x_1' + (\lambda_2-\lambda_1)x_2+ \lambda_3 x_3+\dots+\lambda_n x_n = 0 $$
which is still a convex combination. Note that
$$ 2\lambda_1 f(x_1') + (\lambda_2-\lambda_1)f(x_2) + \lambda_3 f(x_3)+\dots+\lambda_n f(x_n)   $$
$$ \leq \lambda_1(f(x_1)+f(x_2)) + (\lambda_2-\lambda_1)f(x_2) + \lambda_3 f(x_3)+\dots+\lambda_n f(x_n) $$
$$ = \lambda_1 f(x_1) + \lambda_2 f(x_2) + \lambda_3 f(x_3)+\dots+\lambda_n f(x_n)  $$
where we have used $f(x_1') \leq (f(x_1)+f(x_2))/2$ (see the definition of $\varepsilon$-uniform convexity).
The inequality means that $S_1=\{x_1',x_2,\dots,x_n\}$ is an improvement of $S$ in the sense of the approximation to $\breve{f}$. Note also that $x^*(x_1') \leq a/2$.\\ 
In case, $\lambda_1 > \lambda_2$, we will use several vectors $x_k$ with $x^*(x_k) \leq 0$ to compensate $x_1$. This is possible because $a \geq b$ implies that the ``mass'' lying on the halfspace $x^*<0$ is not lesser than $\lambda_1$.
In this case, $\lambda_1$ could be cancelled with several $\lambda_k$'s. In any case, we will get a new set $S_1$ whose cardinal is not larger than that of $S$ and $\conv(S_1) \subset \conv(S)$. After that, suppose that, unfortunately, we still have $\sup x^*(S_1) = a$. In such a case, the maximizing vector cannot be $x_1$, so it is a new vector, say $x_3$. We will apply the argument with $x_3$ in order to replace it by another vector $x_3'$ and $S_1$ by a new set $S_2$. Eventually, we will get $\sup x^*(S_n) \leq a/2$ after a finite number of steps. Then, with the same $x^*$, we have to change the constants $a,b>0$ by new ones. This can be done 
with the same $x^*$ until we get $\max\{a,b\} < \varepsilon$, so it is not possible to go further.\\
If the set of points it is not yet inside $B(0,\varepsilon)$ then find a new $x^* \in S_{X^*}$ such that  $ \sup x^*(S_n) \geq \varepsilon$ and then run again the algorithm.
Since $\conv(S)$ is finite dimensional, it is enough to do this procedure over finitely many $x^* \in S_{X^*}$ in order to get $S_n \subset B(0,\varepsilon)$ eventually.\findemo\\

\noindent
{\bf Proof of Theorem \ref{main-convexify}.}
For the proof it would be convenient to represent a convex combination in $X$ by means of a vector integral instead of the usual symbol ``$\sum$''. Namely, given a convex combination $\sum_{k=1}^n \lambda_k x_k$, the weights $\lambda_k$ are changed by $n$ disjoint intervals $I_k$ of lengths $\lambda_k$ and whose union is $[0,1]$. In this way, the convex combination is represented as the integral of the simple function $\overline{x}$ defined on $[0,1]$ by 
$\overline{x}(t)=x_k$ whenever $t \in I_k$. As we will deal only with simple functions, no further knowledge of vector integration theory is required.\\
We resume the proof. If $\breve{f}$ is proper then it is bounded below by an affine function, so by adding an affine function (that does not alter the $\varepsilon$-uniform convexity), we may suppose that $f$ is bounded below (actually that is true without modifications, see Corollary \ref{coro-proper}).
Given $x,y \in \dom(f)$, if $\|x-y\| \geq 3\varepsilon$ then 
$$\breve{f} \left(\frac{x+y}{2} \right)  \leq \frac{\breve{f}(x)+\breve{f}(y)}{2} - \delta_f(\varepsilon).$$
Indeed, fix $\eta>0$. By Proposition \ref{bounded-reduction}, we may take $(\overline{x}_n)$ a sequence of simple functions defined on $[0,1]$ such that $\| \overline{x}_n(t)-x \| < \varepsilon$ for all $t \in [0,1]$, $n \in {\Bbb N}$, with  
$$ \lim_n \int_0^1 \overline{x}_n(t)\,dt = x  , ~~\mbox{and} ~~ \lim_n \int_0^1 f(\overline{x}_n(t))\,dt = \breve{f}(x).$$ 
Let $(\overline{y}_n)$ and analogous sequence of simple functions playing the same role for $y$ and $\breve{f}(y)$. Clearly we have $\|\overline{x}_n(t)-\overline{y}_n(t)\| \geq \varepsilon$ for all $t \in [0,1]$ and $n \in {\Bbb N}$. 
Therefore 
$$ \breve{f}\left( \frac{x+y}{2} \right) \leq \liminf_n \int_0^1 f \left( \frac{\overline{x}_n(t)+\overline{y}_n(t)}{2} \right) \, dt $$
$$ \leq \lim_n \int_0^1 \left( \frac{f(\overline{x}_n(t))+f(\overline{y}_n(t))}{2}  - \delta_f(\varepsilon)\right) dt \leq \frac{\breve{f}(x) + \breve{f}(y) }{2}
- \delta_f(\varepsilon) .$$
Since $\eta>0$ was arbitrary we get the claimed inequality provided $\|x-y\| \geq 3\varepsilon$.\\
Now we will suppose $\varepsilon \leq \|x-y\| <3\varepsilon$. Proposition \ref{bounded-reduction} implies that reducing the domain of $f$ to $[x,y] + B(0,\varepsilon)$ does not affect to the values of $\breve{f}(x)$, $\breve{f}(y)$ and $\breve{f}(\frac{x+y}{2})$.
Fix $\varepsilon' >\varepsilon$. Theorem \ref{convexify} says that
$\delta_{\breve{f}}(\varepsilon')$ depends only on $\varepsilon'$, $\delta_f(\varepsilon)$, $\varpi$, which are fixed, 
and the diameter of the domain, which is bounded by  $5\varepsilon$.\findemo

\section{Building uniformly convex functions}

Most of the constructions of uniformly convex functions on a Banach spaces that one can find in the literature are based on modifications of a uniformly convex norm, see \cite{BV2}. Nevertheless, the existence of a finite uniformly convex function whose domain has nonempty interior implies that $X$ has an equivalent uniformly convex norm, see \cite[Theorem 2.4]{BGHV}. In any case, the constructions dealing with the composition of a uniformly convex norm and a suitably chosen function can be quite tricky, except for the Hilbert space. Here we will exploit a method based on ``discretization'' and uniformly quasi-convex functions.

\begin{lema}\label{quasi-3}
Let $\varepsilon >0$ and let 
$f: X \rightarrow {\Bbb R}$ be a bounded below $\varepsilon$-uniformly quasi-convex function with modulus $\delta>0$.
Then the function $h \circ f$ is $\varepsilon$-uniformly convex, where
$h(t)=3^{t/\delta}$.
\end{lema}

\noindent
{\bf Proof.}
The function $h$ is increasing and satisfies the property $3h(t) = h(t+\delta)$. Take 
$$ \eta=4^{-1}\inf\{ h(t+\delta) - h(t): t \geq \inf f \} = 2^{-1} \cdot 3^{\inf f /\delta}$$
and note that it depends only on $f$. If $x,y \in \dom(f)$ are such that $d(x,y) \geq \varepsilon$ take
$a=f(x)$, $b=f(y)$ and $c=f(\frac{x+y}{2})$. 
The hypothesis says that $c \leq \max\{a,b\}-\delta$.
With no loss of generality, we may assume $b \leq a$. We have
$$ h(c) \leq h(a) - 4\eta .$$
Since $3h(c) \leq h(a)$ and $h(b)>0$, we also have
$$ 3 h(c) \leq h(a) + 2h(b)  $$
and adding the previous inequality, we get
$$ 4 h(c) \leq 2 h(a) + 2 h(b) -4\eta $$
and thus
$$ h(c) \leq \frac{h(a)+h(b)}{2} - \eta $$
which is the $\varepsilon$-uniform convexity of $h \circ f$.\findemo\\

If $X$ is uniformly convex, it is well known that $x \rightarrow \|x\|^2$ is a uniformly convex function on bounded convex subsets. The usual construction of a global uniformly convex functions involves additional properties of the norm, such as having a power type modulus of uniform convexity. Here there is a simple alternative construction based in our methods.

\begin{prop}
If $X$ has a uniformly convex norm then there exists a real function $\phi$ such that $x \rightarrow \phi(\|x\|)$ is a uniformly convex function defined on $X$.
\end{prop}

\noindent
{\bf Proof.}
Fix $\varepsilon>0$. Take $a_1=\varepsilon/2$ and define inductively a sequence $(a_n)$ by the implicit equation
$$ a_{n-1} = \left( 1 - \delta_X  \left(  \frac{\varepsilon}{a_{n}} \right) \right) a_{n} $$
which has a unique solution thanks to the continuity of $\delta_X$ on $[0,2)$, \cite[Lemma 5.1]{GoKi}. The sequence $(a_n)$ is increasing with 
$\lim_n a_n=+\infty$ and has the following property:  if $\|x\|,\|y\| \leq a_n$ and $\|x-y\| > \varepsilon$ then $\| (x+y)/2 \| \leq a_{n-1}$.\\
Define a function  as $f_\varepsilon(x)=n$ if $a_{n-1} < \|x\| \leq a_n$. Note that $f_\varepsilon$ satisfies the hypothesis of Lemma \ref{quasi-3} with $\delta=1$, and so $h \circ f_\varepsilon$ is $\varepsilon^+$-uniformly convex. Now, for $\varepsilon=1/n$, take
$f_n$ the convex envelope of $h \circ f_\varepsilon$ and $c_n =2^{-n} (\sup f_n( nB_X))^{-1}$. The series $\sum_{n=1}^\infty c_n f_n$ converges uniformly on bounded sets to a uniformly convex function $f$. 
By construction, $f(x)$ depends only on $\|x\|$. Therefore, we may define a real function by $\phi(t)=f(x)$ if $t=\|x\|$, for $t \geq 0$. Clearly, $f(x)=\phi(\|x\|)$.\findemo\\

Now we will explain constructions using trees. The definition of $\varepsilon$-separated (dyadic) tree was given in the introduction. Bushes are defined in a very similar way, however the index set is $\bigcup_{k=0}^n {\Bbb N}^k$ and
$x_s=\sum_k \lambda_{s\frown k}  x_{s\frown k} $ where $ \lambda_{s\frown k} \geq 0$, $ \lambda_{s\frown k} =0$ except for finitely many $k$'s and $ \sum_k \lambda_{s\frown k} =1$. We say that a bush $\{ x_s: |s| \leq n \}$ is 
$\varepsilon$-separated if $\| x_{s\frown k} - x_{s} \| \geq \varepsilon $ for all $k$ such that  $\lambda_{s\frown k}>0$.
In this way, an $\varepsilon$-separated tree is a particular case of an $\varepsilon/2$-separated bush.
Separated trees and bushes are obstructions to the existence of bounded uniformly convex functions. 

\begin{prop}\label{bush}
Let $\varepsilon>0$ and let $C \subset X$ be a convex set that supports an  $\varepsilon$-uniformly convex function $f$ with values in 
$[a,b]$. Then
$(b-a)/\delta_f(\varepsilon)$ is the maximum height of
\begin{enumerate}
\item any $\varepsilon$-separated tree contained in $C$;
\item any $\varepsilon^+$-separated bush contained in $C$.
\end{enumerate}
\end{prop}

\noindent
{\bf Proof.} If $\{ x_s\}$ is an  $\varepsilon$-separated  tree then we have
$$f(x_s) \leq \max\{ f(x_{s\frown 0}),  f(x_{s\frown 1})  \}-\delta_f(\varepsilon)$$
that gives the estimation. In the case of bushes, the argument is the same after passing to $\breve{f}$, which is $\varepsilon^+$-uniformly convex by
Theorem \ref{main-convexify}, and applying Corollary \ref{convex-delta}.\findemo\\

Our following result is quite a converse.

\begin{theo}\label{tree}
Let $\varepsilon>0$ and let $C \subset X$ be a  convex set such that contains not arbitrarily high $\varepsilon$-separated  trees (with respect some uniformly continous pseudometric). Then $C$ supports a bounded $\varepsilon$-uniformly convex function, and a bounded convex $\varepsilon^+$-uniformly convex function (with respect the same pseudometric).
\end{theo}

\noindent
{\bf Proof.}
Define a function for $x \in C$ by
$$ f(x) = - \max\{ \mbox{height}(x_s): (x_s) \subset C\, \, \varepsilon \mbox{-sep. tree}, x_\emptyset =x \} $$
and $f(x)=+\infty$ otherwise. We claim that $f$ 
is $\varepsilon$-uniformly quasi-convex.
Indeed, consider points $x,y \in C$ with $d(x,y) \geq \varepsilon$. Take $\varepsilon$-separated trees contained into $C$
$\{ x_{s'}: |s'| \leq n' \}$ and $\{ y_{s''}: |s''| \leq n'' \}$ 
of maximal length with the property that $x_\emptyset=x$ and $y_\emptyset=y$.
The trees can be ``glued'' as follows. Take $n=\min\{n',n''\}$. Define a new tree $(z_s)$, for $|s| \leq n+1$, as $z_\emptyset=\frac{x+y}{2}$, $z_{0\frown s}= x_s$ and $z_{1\frown s}= y_s$ for $|s| \leq n$.
Now $(z_s)$ is a $\varepsilon$-separated tree rooted at $\frac{x+y}{2}$ of height $\min\{n',n''\}+1$. 
That means in terms of the function $f$ the uniform quasi-convex inequality
$$ f \left( \frac{x+y}{2} \right) \leq \max\{ f(x), f(y)\} -1 $$
for $d(x,y) \geq \varepsilon$. Now, 
Lemma \ref{quasi-3} says that $h \circ f$ is $\varepsilon$-uniformly convex and its convex hull is $\varepsilon^+$-uniformly convex after  Theorem \ref{main-convexify}.\findemo\\

\noindent
{\bf Proof of Theorem \ref{main-swc}.}  It just follows from Theorem \ref{tree} and Proposition \ref{bush}.\findemo\\

Finally we will explain constructions based on the dentability index.
Let $C$ be a bounded closed convex set of $X$, $(M,d)$ a pseudometric space and $F: C \rightarrow M$ a map. 
We say that $F$ is {\it dentable} if for any nonempty closed convex subset $D \subset C$ and
$\varepsilon >0$, it is possible to find an open halfspace $H$ intersecting $D$ such that
$\mbox{diam}(F(D \cap H)) <\varepsilon$, where the notation ``$\diam$'' stands for the diameter is computed with respect to $d$. 
If $F$ is dentable, we may consider
the following {\it set derivation}
\[ [D]'_{\varepsilon} = \{ x \in D: \mbox{diam}(F(D \cap H))>\varepsilon, ~ \forall \, H \in {\Bbb H}, \,
x \in H \} .\]
Here ${\Bbb H}$ denotes the set of all the open halfspaces of $X$. Clearly, $[D]'_{\varepsilon}$
is what remains of $D$ after removing all the slices whose diameter through $F$ is less or equal than
$\varepsilon$. A useful trick is the so called (nonlinear) Lancien's midpoint argument: if a segment 
satisfies $[x,y] \subset D$ and $[x,y] \cap [D]'_{\varepsilon} = \emptyset$ then
$d(F(x),F(y)) \leq 2\varepsilon$, see the beginning of \cite[Theorem 2.2]{raja}.
Consider the sequence of sets defined by $[C]_{\varepsilon}^{0}=C$ and, for every $n \in {\Bbb N}$, inductively by
\[ [C]_{\varepsilon}^{n}=[[C]_{\varepsilon}^{n-1}]'_{\varepsilon}. \]
If there is $n \in {\Bbb N}$ such that $[C]_{\varepsilon}^{n-1} \not = \emptyset$ and  $[C]_{\varepsilon}^{n}  = \emptyset$ we say that $\Dz(F,\varepsilon)=n$. We say that $F$ is {\it finitely dentable} if $\Dz(F,\varepsilon) <  \omega$ for every 
$\varepsilon>0$ ($ \omega$ stands for the first infinite ordinal number).
All these notions can be applied to the identity map of a convex set where there is a pseudometric defined. The following result is the quantified version of \cite[Theorem 2.2]{raja}. For convenience we will write
$$ \Delta_\Phi(x,y) =  \frac{\Phi(x)+\Phi(y)}{2} - \Phi \left( \frac{x+y}{2} \right) . $$

\begin{theo}\label{delta-renor}
Let $C \subset X$ be a bounded closed convex set, let $M$ be a pseudometric space, let
$F:C \rightarrow M$ be a uniformly continuous map, and  let $\varepsilon>0$.
\begin{enumerate}
\item Suppose that there exists a bounded lower semi-continuous convex function $\Phi$ defined on $C$ and $\delta>0$ such that
$d(F(x),F(y)) \leq \varepsilon$ whenever $x,y \in C$ satisfy  $\Delta_\Phi(x,y) \leq \delta$. Then
$\Dz(F,\varepsilon)< \omega$. 
\item On the other hand, if $\Dz(F,\varepsilon)< \omega$ then for every $\varepsilon'>2\varepsilon$
there exits a bounded lower semi-continuous convex function 
$\Phi$ defined on $C$ and $\delta>0$ such that $d(F(x),F(y)) \leq \varepsilon'$ whenever
$x,y \in C$ satisfy $\Delta_\Phi(x,y) \leq \delta$.
\end{enumerate}
\end{theo}

\noindent
{\bf Proof.}
Let $s=\sup f(C)$. The hypothesis implies $ [C]_{\varepsilon}' \subset \{f \leq s-\delta\}$. Iterating this we will eventually get to the empty set. For the second part, we need to introduce some notation. 
Firstly put $d'(x,y)=d(F(x),F(y))$ which is a pseudometric uniformly continuous with respect to $\| \cdot \|$. Derivations and diameters will be referred to $d'$.
The slice of a set $A$ with parameters $x^* \in X^*$ and 
$\alpha>0$ is
$$ S(A,x^*,\alpha) = \{ x \in A: x^*(x) > \sup x^*(A) -\alpha \} .$$
The ``half-derivation'' of a convex set is defined as 
\[ \langle D \rangle'_{\varepsilon} = \{ x \in D: x^*(x) \leq \alpha, ~ \forall x^*,\alpha>0 ~\mbox{such that}~  
\diam(S(D,x^*,2\alpha)) > \varepsilon \} .\]
The geometric interpretation is that we remove half of the slice, in sense of the width, for every slice of $d'$-diameter less than $\varepsilon$. This derivation can be iterated by taking $\langle C \rangle_{\varepsilon}^{n}=\langle \langle C \rangle_{\varepsilon}^{n-1} \rangle'_{\varepsilon}$.
It is not difficult, but rather tedious, to show that if $\Dz(F,\varepsilon)< \omega$ then for every $\varepsilon'>2\varepsilon$ there is some $N \in {\Bbb N}$ such that $\langle C \rangle_{\varepsilon'}^{n}=\emptyset$. The idea is the following. Firstly note that every slice of $C$ not meeting $[C]'_{\varepsilon}$ has diameter $2\varepsilon$ at most by Lancien's argument.
Taking ``half a slice'' of the slice given by some $x^* \in X^*$, we deduce that 
$$ \sup x^*( \langle C \rangle'_{2\varepsilon} ) - \sup x^*( [C]'_{\varepsilon} ) \leq 
2^{-1} ( \sup x^*(  C ) - \sup x^*( [C]'_{\varepsilon} )). $$
Iterating, we would get
$$ \sup x^*(  \langle C \rangle^n_{2\varepsilon}  ) - \sup  x^*( [C]'_{\varepsilon} ) \leq 
2^{-n} ( \sup x^*(  C ) - \sup  x^*( [C]'_{\varepsilon} ))$$
for every $x^* \in X^*$.
If $\eta>0$, we will get for some $n$ large enough  that
$$  \langle C \rangle^n_{2\varepsilon}  \subset [C]'_{\varepsilon} +B(0,\eta). $$
We can do that for every set $[C]^k_{\varepsilon} $.  A perturbation argument, using
 the room between $\varepsilon$ and $\varepsilon'$, will 
allow us to fill the gap between the sequences of sets. In this way  we will get that $ \langle C \rangle^n_{\varepsilon'}=\emptyset$ for some $n \in {\Bbb N}$ large enough.\\
Now we define a function $g$ on $C$ by 
$g(x)=-n$ if $x \in \langle C \rangle_{\varepsilon'}^{n} \setminus \langle C \rangle_{\varepsilon'}^{n+1}$ 
following the notation above.
We claim that $g$ satisfies Lemma \ref{quasi-3} with separation $\varepsilon'$. 
Indeed, if $d'(x,y) > \varepsilon'$ and $n=-\max\{g(x),g(y)\}$ then $x,y \in \langle C \rangle_{\varepsilon'}^{n}$. 
If $\frac{x+y}{2} \not \in \langle C \rangle_{\varepsilon'}^{n+1}$ then the segment $[x,y]$ would be fully contained into a slice of diameter less than $\varepsilon'$ and so $d'(x,y) \leq \varepsilon'$ which is a contradiction.
Therefore $\frac{x+y}{2} \in \langle C \rangle_{\varepsilon'}^{n+1}$ and so $g(\frac{x+y}{2}) \leq -n-1$. Now  $f(x)=3^{g(x)}$ is $\varepsilon'$-uniformly convex with respect to $d'$. Take $\Phi=\breve{f}$ to get the desired function.\findemo\\

If $F$ in Theorem \ref{delta-renor} (2) were finitely dentable, a standard argument using a  convergent  series would lead to this results, which is essentially \cite[Proposition 3.2 ]{GL} with a uniformly convex function instead of a norm.

\begin{coro}
Let $C \subset X$ be a bounded closed convex set, let $M$ be a pseudometric space, and let
$F:C \rightarrow M$ be a uniformly continuous finitely dentable map.
Then there exits a bounded convex
function $\Phi$ defined on $C$ such that for every $\varepsilon>0$ there is $\delta>0$ such that $d(F(x),F(y)) \leq \varepsilon$ whenever $x,y \in C$ are such that $\Delta_\Phi(x,y) \leq \delta$.

\end{coro}

\section{Improving functions and domains}

So far, the best improvement we have done on an existing $\varepsilon$-uniformly convex function is taking its lower semicontinuous convex envelope provided this last one is proper. The aim in this section is to manipulate the functions in order to improve their qualities. We will begin by proving the results about global behaviour.\\

\noindent
{\bf Proof of Theorem \ref{main-coercive}}.
Since $\breve{f} \leq f$, it is enough to prove that the property holds for an $\varepsilon$-convex and convex proper function. Actually the same proof for a uniformly convex function done in Zalinescu's book \cite[Proposition 3.5.8]{Za2} works in this case because $\liminf_{t\rightarrow +\infty} t^{-2} p_{f}(t) \geq \varepsilon^{-2} p_{f}(\varepsilon)>0$.\\
For the second part,
without loss of generality we may assume that $x_0^*=0$ (just change $f$ by $f+x_0^*$). 
Let $\delta=\delta_{\breve{f}}(\varepsilon')$ and take $\eta=\inf f+\delta-\breve{f}(x_0)>0$.
Note that $\inf f=\inf \breve{f} $. By the property established in the first part applied to $\breve{f}-x^*$, there is $R>0$ such that $\breve{f}(x) \geq \breve{f}(x_0)-x^*(x-x_0)$ for any $x^* \in B_{X^*}$ and $\|x-x_0\| \geq R$. 
Now, fix $x^*$ such that $\|x^*\| \leq \eta/R$. Then we have
$$ \breve{f}(x) +x^*(x) \geq  \breve{f}(x_0) + x^*(x_0) -\delta $$
for all $x \in X$ such that $\|x-x_0\| \leq R$, and therefore the inequality holds for all $x \in X$. 
That implies $\epi(\breve{f} +x^*+\delta)$ does not meet the horizontal slice 
$$ S =\{(x,t) \in \epi(\breve{f} + x^*): t \leq \breve{f}(x_0)+x^*(x_0)\} $$
By Proposition \ref{criterio-slice}, the projection of $S$ on $X$ has diameter less than $\varepsilon'$.
Moreover, if $f+x^*$ attains a minimum at $x$, then the same holds for $\breve{f}+x^*$ and so $x \in S$. Since $x_0 \in S$ we have
$\|x-x_0\| \leq \varepsilon'$. The existence of a dense set of $x^*$'s such that $\breve{f}+x^*$ attains a minimum is guaranteed by Brøndsted--Rockafellar \cite[Theorem 4.3.2]{BV} (or Bishop--Phelps \cite[Theorem 7.4.1]{banach} applied to the epigraph).\findemo\\

As a consequence, we characterize when an $\varepsilon$-uniformly convex function has a proper convex envelope.

\begin{coro}\label{coro-proper}
Let $\varepsilon>0$ and let $f$ be an $\varepsilon$-uniformly convex function. Then the following statements are equivalent:
\begin{enumerate}
\item $\breve{f}$ is proper;
\item $f$ is bounded below;
\item  $f$ is bounded below by an affine continuous function.
\end{enumerate}
\end{coro}

For a  $\varepsilon$-uniformly quasi-convex function we can say the following

\begin{prop}
Let $\varepsilon>0$ and let $f$ be an $\varepsilon$-uniformly quasi-convex function that is bounded below. 
Then $f$ is coercive and moreover
$$ \liminf_{\|x\| \rightarrow +\infty} \frac{f(x)}{ \|x\|} >0. $$
\end{prop}

\noindent
{\bf Proof.} By adding a constant, we may suppose that $\inf f=0$. Take $x_0 \in X$ such that $f(x_0) < \delta/2$. For any $x \in X$ such that $\|x-x_0\| \geq \varepsilon$ we have $f(x) \geq \delta$. Indeed, otherwise it would be $f(x) < \delta$ and by the $\varepsilon$-uniformly quasi-convexity, $f(\frac{x+x_0}{2}) < \inf f$, an obvious contradiction.
Now, if $\|x-x_0\| \geq 2\varepsilon$, then $\| \frac{x+x_0}{2} - x_0 \| \geq \varepsilon$. That implies $f( \frac{x+x_0}{2}) \geq \delta$ and therefore $f(x) \geq 2\delta$. Inductively, we will get that if $\|x-x_0\| \geq 2^n \varepsilon$ then $f(x) \geq 2^n \delta$. Now, the statement follows easily.\findemo\\

The following results will show that, given a
$\varepsilon$-uniformly convex function, we can make modifications in both the function and its domain in order to get a new function with additional properties.

\begin{prop}
Let $\varepsilon>0$, let $f$ be an $\varepsilon$-uniformly convex function that is locally bounded below and let $\eta>0$. Then there exists a lower semicontinuous $(\varepsilon+2\eta)$-uniformly convex function defined on $\dom(f) + B(0,\eta)$. In particular, $\overline{ \dom(f) }$ admits a lower semicontinuous $\varepsilon^+$-uniformly convex function. 
\end{prop}

\noindent
{\bf Proof.} Define $g(x)= \inf\{f(y): \|y-x\| < \eta\}$ on $\dom(f) + B(0,\eta)$. 
This function $g$ is $(\varepsilon+2\eta)$-uniformly convex (the simple verification of this fact is left to the reader).
Now take its lower semicontinuous envelope.\findemo\\

The following result will be done for  $\varepsilon$-uniformly convexity with respect to a metric because such a degree of generality will be needed later.

\begin{prop}\label{add-Lipschitz}
Let $\varepsilon>0$, let $f$  be an $\varepsilon$-uniformly convex function (with respect to a metric $d$ with modulus of uniform continuity $\varpi$) and  let $C \subset \dom(f) $ be convex such that $f$ is bounded on it. Then for any $\varepsilon'>\varepsilon$, there exists $g \in \Gamma(X)$ Lipschitz (with respect to the norm of $X$) such that $g|_C$ is $\varepsilon'$-uniformly convex.
\end{prop}

\noindent
{\bf Proof.} Without loss of generality, we may assume that $f$ is convex.
Take $\eta>0$ such that $\varpi(\eta)<(\varepsilon'-\varepsilon)/2$,
$m=\sup\{f(x)-f(y): x,y \in C\}$ and $c=m/\eta$. Define
$$ g(x) = \inf\{ f(y)+c\|x-y\|: y \in C \} $$
which is convex and $c$-Lipschitz. Let $x \in C$ and $\xi>0$. Then either $g(x)=f(x)$ and the infimum is attained with $y=x$, or $g(x) < f(x)$. In the last case, the infimum can be computed over the $y \in C$ such that $f(y)+c\|x-y\| < f(x)$. 
Therefore, we can find $y \in C$ such that $ f(y)+c\|x-y\| < g(x) + \xi$ and $\|x-y\| < m/c=\eta$, which implies
$d(x,y)<(\varepsilon'-\varepsilon)/2$.
Now, for $x_1,x_2 \in C$ with $d(x_1,x_2) \geq \varepsilon'$ find $y_1,y_2 \in C$ as above.
Clearly we have $d(y_1,y_2) \geq \varepsilon$, and so
$$ g \left( \frac{x_1+x_2}{2} \right) \leq f \left( \frac{y_1+y_2}{2} \right) + c \left\| \frac{x_1+x_2}{2} - 
\frac{y_1+y_2}{2} \right\| $$
$$ \leq \frac{f(y_1)+f(y_2)}{2} - \delta + \frac{c}{2}  \| x_1-y_1 \| +   \frac{c}{2}   \|x_2-y_2 \|   \leq  \frac{g(x_1)+g(x_2)}{2} - \delta +2\xi .$$
Since $\xi>0$ was arbitrary, we get the $\varepsilon'$-uniform convexity of $g$ as wished.\findemo

\begin{rema}
A Baire category argument shows that an $\varepsilon$-uniformy convex function $f$ is bounded in an open ball if $\dom(f)$ has nonempty interior. However, we do not know how to ensure that $f$ will be bounded on a larger set.
\end{rema}

Now we will show how to change an $\varepsilon$-uniformly convex function by a norm with the same property.

\begin{theo}\label{renorming}
Let $(X,\| \cdot \|)$ be a Banach space,
let $f \in \Gamma(X)$ be a non negative function and let $C \subset \dom(f)$ be a bounded convex set.
Assume $f$ is Lipschitz on $C$.
Then given $\delta>0$ there exists an equivalent norm $|\!|\!|  \cdot |\!|\!| $ on $X$ and $\zeta>0$ such that 
$\Delta_f(x,y) <\delta$ whenever $x,y \in C$ satisfy $\Delta_{ |\!|\!|  \cdot |\!|\!|^2 }(x,y) <\zeta$.
Therefore, if $f$ 
was moreover $\varepsilon$-uniformly convex for some $\varepsilon>0$ (with respect to a pseudometric) on $C$, then $|\!|\!|  \cdot |\!|\!|^2$ would be 
$\varepsilon$-uniformly convex on $C$ (with respect to the same pseudometric).
\end{theo}

\noindent
{\bf Proof.}
Taking $f(x)+f(-x)+\|x\|$ instead, we may indeed assume that $f$ is symmetric and attains a strong minimum at $0$. 
Let $M=\sup f(C)$ and $m=\min f(C)+\delta/2$. The Lipschitz condition easily implies that there is $\eta>0$ such that if $r \leq M$ then 
$$\{f \leq r-\delta\} + B(0,\eta) \subset \{f \leq r\}.$$
For $r \in [m,M]$ let $\| \cdot \|_r$ the Minkowski functional of the set $\{f \leq r\}$, which is an equivalent norm on $X$. 
Let $N=\sup\{ \|x\|:x \in C\}$ and note that $\lambda=(1+\eta/N)^{-1}$ has the property
$$ \{f \leq r-\delta\}  \subset \lambda \, \{f \leq r\} .$$
We deduce the following property:
if $x,y \in C$, $\|x\|_r,\|y\|_r \leq 1$ and $\Delta_f(x,y) \geq \delta$
 then 
$$ \left\| \frac{x+y}{2} \right\|_r \leq \lambda. $$
Consider a partition $m=a_1<a_2<\dots<a_k=M$ such that $a_{j} /a_{j+1}<\lambda^{1/2}$ and put $|\!|\!|  \cdot |\!|\!| _j=\| \cdot \|_{a_j}$. Let $x,y \in C$ such that $\Delta_f(x,y) \geq \delta$.
Assume $f(x) \geq f(y)$ for instance. 
There is some $1 \leq j<k$ such that $1 \geq |\!|\!| x|\!|\!| _j \geq \lambda^{1/2}$.
Since $ |\!|\!| \frac{x+y}{2}  |\!|\!|_j \leq \lambda$, we have
$$  \left|\!\left|\!\left| \frac{x+y}{2}  \right|\!\right|\!\right|_j \leq 
\max\{ |\!|\!| x|\!|\!| _j , |\!|\!| y |\!|\!| _j  \} - (\lambda^{1/2} -\lambda). $$
Following the same computations that in  the proof of Proposition \ref{square-quasi}, we have
$$  \left|\!\left|\!\left| \frac{x+y}{2}  \right|\!\right|\!\right|^2_j \leq 
\frac{ |\!|\!| x|\!|\!|_j^2 +  |\!|\!| y |\!|\!|^2_j  }{2} - \frac{(\lambda^{1/2} -\lambda)^2}{4}. $$
Therefore, if we define an equivalent norm by 
$|\!|\!| \cdot |\!|\!|^2 =\sum_{j=1}^k |\!|\!|  \cdot |\!|\!|^2 _j$
we will have
$$  \left|\!\left|\!\left| \frac{x+y}{2}  \right|\!\right|\!\right|^2 \leq 
\frac{ |\!|\!| x|\!|\!|^2 +  |\!|\!| y |\!|\!|^2  }{2} - \frac{(\lambda^{1/2} -\lambda)^2}{4}. $$
whenever $x,y \in C$ satisfies $\Delta_f(x,y) \geq \delta$, meaning that the statement is true with
$\zeta=4^{-1} (\lambda^{1/2} -\lambda)^2$.\findemo\\

\noindent
{\bf Proof of Theorem \ref{main-renorming}.}  
Consider the sets $C_n=\{f \leq n\}$ that eventually will ``capture'' any set where $f$ is bounded. Fixed $n \in {\Bbb N}$,
by  Proposition \ref{add-Lipschitz}, we may assume that $f$ is already norm-Lipschitz and finite on $X$ provided we change $\varepsilon$ by $\varepsilon^+$.
Let $|\!|\!| \cdot |\!|\!|_n$ the norm given by Theorem \ref{renorming}, which is $\varepsilon^+$-uniformly convex on $C_n$. 
Let $\alpha>0$ and let 
$(\alpha_n)$ be a sequence of positive numbers such that
$$|\!|\!| \cdot |\!|\!|^2 = \alpha |\!|\!| \cdot |\!|\!|^2 + \sum_{n=1}^\infty \alpha_n |\!|\!|  \cdot |\!|\!|^2 _n$$
converges uniformly on bounded sets. Clearly, $|\!|\!| \cdot |\!|\!|$ will be $\varepsilon^+$-uniformly convex on bounded sets too. 
The last affirmation follows just taking $\alpha>0$ small enough.\findemo\\

Finally we will discuss the approximation by differences of convex functions in terms of the index of dentability improving \cite[Theorem 1.4]{raja} and \cite[Theorem 4.1]{GL}, which in turn are based in the seminal work \cite{cepe}.
A real function defined on a convex set is called ${\mathcal DC}$-Lipschitz if it is the difference of two convex Lipschitz functions with the same domain. For convenience, we have to make explicit the domain $C$.
The symbol $\| \cdot \|_C$ stands for the supremum norm on $C$.

\begin{theo}
Let $C \subset X$ be a bounded closed convex set and let $f : C \rightarrow {\Bbb R}$ be a uniformly continuous function. Consider the following numbers:
\begin{itemize}
\item[$(\varepsilon_1)$] the infimum of the $\varepsilon>0$ such that $\Dz(f,\varepsilon) <\omega$;
\item[$(\varepsilon_2)$] the  infimum of the $\varepsilon>0$ such that there exists a ${\mathcal DC}$-Lipschitz function $g$ such that $\|f-g\|_C < \varepsilon$.
\end{itemize}
Then   $\varepsilon_1/2 \leq \varepsilon_2 \leq 2\varepsilon_1$.
\end{theo}

\noindent
{\bf Proof.} 
Let $\varepsilon>\varepsilon_2$ and find a ${\mathcal DC}$-Lipschitz function $g$ such that $\|f-g\|_C < \varepsilon$. 
We know by \cite[Proposition 5.1]{raja} that $g$ is finitely dentable, which easily implies that $f$ is $2\varepsilon$-finitely dentable.\\
For the reverse inequality, 
take $\varepsilon>2\varepsilon_1$ and $M=\sup\{f(x)-f(y):x,y \in C\} <+\infty$.
Apply Theorem \ref{delta-renor} to get a function $\Phi$ such that $|f(x)-f(y)| \leq \varepsilon$ if $\Delta_\Phi(x,y) <\delta$.
By Proposition \ref{add-Lipschitz} we may suppose that $\Phi$ is Lipschitz too, and by Theorem \ref{renorming}, there is an  equivalent norm $ |\!|\!| \cdot  |\!|\!|$ defined on $X$ such that 
$\Delta_{  |\!|\!| \cdot  |\!|\!|^2 }(x,y) <\zeta$ implies $\Delta_\Phi(x,y) <\delta$.
Take $c > M/\zeta$. Consider the function 
$$ g(x) = \inf_{y \in C} \left\{ f(y) + c\left( \frac{ |\!|\!| x  |\!|\!|^2 + |\!|\!| y  |\!|\!|^2}{2}- 
\left|\!\left|\!\left| \frac{x+y}{2}  \right|\!\right|\!\right|^2 \right) \right\} = 
 \inf_{y \in C}  \{  f(y) + c \, \Delta_{  |\!|\!| \cdot  |\!|\!|^2 }(x,y) \} $$
which is actually an $\inf$-convolution with the {\it Cepedello's kernel}, see \cite{cepe} or \cite[Theorem 4.21]{BL}.
For every $x \in C$, the infimum can be computed just on the set
$$ A(x) = \{ y \in C:  f(y) + c \, \Delta_{  |\!|\!| \cdot  |\!|\!|^2 }(x,y) \}  \leq f(x) \} .$$
If $x \in C$ and $y \in A(x)$, we have
$$  0 \leq  c \, \Delta_{  |\!|\!| \cdot  |\!|\!|^2 }(x,y)  \leq f(x)-f(y) \leq M .$$
Then 
$  \Delta_{  |\!|\!| \cdot  |\!|\!|^2 }(x,y)  \leq \zeta$ by the choice of $c$ and thus $0 \leq f(x)-f(y) \leq \varepsilon$. 
Fix $\eta>0$ and take $y \in A(x)$ such that
$$  f(y) +  c \, \Delta_{  |\!|\!| \cdot  |\!|\!|^2 }(x,y)  \leq g(x) + \eta . $$
Then
$$ f(x) - g(x) \leq f(x) - f(y) - c \, \Delta_{  |\!|\!| \cdot  |\!|\!|^2 }(x,y)  + \eta \leq \varepsilon + \eta. $$
We deduce that $\|f(x)-g(x)\|_C \leq \varepsilon $ and
$$  g(x) = \frac{c}{2} \, |\!|\!|x |\!|\!|^2 - \sup_{y \in C} \left\{ c   \left|\!\left|\!\left| \frac{x+y}{2}  \right|\!\right|\!\right|^2
-  \frac{c}{2} \, |\!|\!| y  |\!|\!|^2 -f(y) \right\}  $$
which is an explicit decomposition of $g$ as a difference of two convex Lipschitz functions on $C$, as wanted.\findemo

\section{Quantifying the super weakly compactness}

The notion of super weak compactness was introduced in \cite{Cheng}. However, some results were established independently in \cite{raja} for an equivalent notion (finite dentable) within the convex sets. Here we will use a definition based on ultrapowers. Given a free ultrafilter ${\mathcal U}$ on a set $\Bbb N$, recall that $X^{\mathcal U}$ is the quotient of $\ell_\infty(X)$ by the subspace of those $(x_n)_{n \in {\Bbb N}}$ such that $\lim_{n, \mathcal U} \|x_n\|=0$. We take $K^{\mathcal U}$ as the image of $K$ by the canonical embedding $x \mapsto  (x)_{n \in {\Bbb N}}$.
A subset $K \subset X$ is said to be {\it relatively super weakly compact} if $K^{\mathcal U}$ is a relatively weakly compact subset of $X^{\mathcal U}$ for a (equivalently all) free ultrafilter ${\mathcal U}$ on ${\Bbb N}$, and $K$ is said to be {\it super weakly compact} if it is moreover weakly closed. The following result gathers several equivalent properties, see \cite{Cheng, raja, raja2, LR}, in order to compare with their quantified versions (Theorem \ref{quantify-swc}).

\begin{theo}\label{known-swc}
Let $C \subset X$ be a bounded closed convex subset. The following statements are equivalent:
\begin{itemize}
\item[$(1)$]  Given $\varepsilon>0$ it is not possible to find arbitrarily long sequences $x_1,\dots,x_n \in C$ such that
$$\dd( \conv\{x_1,\dots,x_k\},  \conv\{x_{k+1},\dots,x_n\}) \geq \varepsilon$$ 
for all $k=1,\dots, n-1$ ($\dd$ stands for the norm distance between two sets);
\item[$(2)$] $C$ contains not arbitrarily high $\varepsilon$-separated dyadic trees for every $\varepsilon>0$;
\item[$(3)$]  $C^{\mathcal U}$ is relatively weakly compact in $X^{\mathcal U}$, for ${\mathcal U}$ a free ultrafilter 
on ${\Bbb N}$ (equivalently, $C$ is relatively SWC, by definition);
\item[$(4)$]  $C^{\mathcal U}$ is dentable in $X^{\mathcal U}$, for ${\mathcal U}$ a free ultrafilter on ${\Bbb N}$;
\item[$(5)$]  $C$ is  finitely dentable;
\item[$(6)$]  $C$ supports a convex bounded uniformly convex function.
\end{itemize}
\end{theo}

In order to state our results we need to introduce some quantities related to sets in Banach spaces. Firstly, a measure of non weakly compactness that has been studied in several papers \cite{FHMZ, granero, CMR}, see also \cite[Section 3.6]{PVJV}. Given $A\subset X$ consider $A \subset X^{**}$ by the natural embedding and take
$$ \gamma(A) = \inf\{ \varepsilon>0: \overline{A}^{w^*} \subset X + \varepsilon B_{X^{**}}\} .$$
It turns out that $\gamma(A)=0$ if and only if $A$ is relatively weakly compact, thus $\gamma$ quantifies the non-weakly compactness of subsets in $X$. This measure is considered more suitable than De Blasi's measure for some problems in Banach space theory. Given a convex set $A \subset X$, let us denote by $\mbox{Dent}(A)$ the infimum of the numbers $\varepsilon>0$ such that $A$ has nonempty slices contained in balls of radius less than $\varepsilon$, and take 
$\Delta(A)= \sup\{ \mbox{Dent}(C): C \subset A\}$. The measure $\mbox{Dent}$ was introduced in \cite{CPR} in relation with the quantification of the RNP property, and actually we may think of $\Delta$ as a measure of non RNP.
Note that $(X^*)^{\mathcal U}$ for  ${\mathcal U}$ a free ultrafilter on ${\Bbb N}$ can be identified with an $1$-norming subspace of $(X^{\mathcal U})^*$ by means of
$$ \langle (x_n^*), (x_n) \rangle = \lim_{n,{\mathcal U}} x_n^*(x_n), $$
where the notation $\langle , \rangle$ is intended to avoid the confusion of dealing with too many brackets.

\begin{lema}\label{ultra-eps}
Let $A \subset X$ be a closed convex bounded subset, ${\mathcal U}$ a free ultrafilter on ${\Bbb N}$ and $\varepsilon>0$. Then
$$  [A^{\mathcal U}]'_{2\varepsilon}   \subset  ([A]'_\varepsilon)^{\mathcal U}   .$$
\end{lema}

\noindent
{\bf Proof.}  Given $(x_n) \in A^{\mathcal U} \setminus  ([A]'_\varepsilon)^{\mathcal U} $, we have to find a slice of 
$A^{\mathcal U}$ containing $(x_n)$ of diameter not greater than $2\varepsilon$. As 
$(x_n) \not \in ([A]'_\varepsilon)^{\mathcal U} $, then for some $\alpha>0$
$$ \{ n : d(x_n, [A]'_\varepsilon) \geq \alpha\} \in {\mathcal U} .$$
Indeed, otherwise the sequence $(x_n)$ would be equivalent to a sequence in $ [A]'_\varepsilon$.
It is possible to find $x_n \in B_{X^*}$ such that $x^*_n(x_n) \geq \alpha + \sup x^*_n ( [A]'_\varepsilon )$ for those indices $n$ from the previous set, for the other $n$'s the choice of  $x_n \in B_{X^*}$ does not make a difference. Consider the functional $(x^*_n) \in (X^*)^{\mathcal U} \subset (X^{\mathcal U})^*$. By construction, 
$$\langle (x^*_n), (x_n) \rangle \geq \alpha +
\sup\langle (x^*_n), ( [A]'_\varepsilon  )^{\mathcal U} \rangle.$$
Now, we will estimate the diameter of the slice defined by $(x^*_n)$.
Suppose that $(y_n), (z_n) \in A^{\mathcal U}$ and
$$ \min\{ \langle (x_n^*), (y_n) \rangle ,  \langle (x_n^*), (z_n) \rangle \} 
\geq \alpha + \sup\langle (x^*_n), ( [A]'_\varepsilon  )^{\mathcal U}  \rangle.  $$
Then, for a subset in $\mathcal U$ of indices $n$, we have $y_n,z_n \in A \cap \{x: x_n^*(x) \geq \alpha + \sup x^*_n (A_n)\}$ and thus $\|y_n-z_n\| \leq 2\varepsilon$ by Lancien's midpoint argument. That implies $\| (y_n) - (z_n) \| \leq 2\varepsilon$, so the diameter of the slice is not greater than $2\varepsilon$ as wished.\findemo\\

The following result is the quantitative counterpart of Theorem \ref{known-swc}.

\begin{theo}\label{quantify-swc}
Let $C \subset X$ be a bounded closed convex subset. Consider the following numbers:
\begin{itemize}
\item[$(\mu_1)$] the supremum of the numbers $\varepsilon>0$ such that for any $n \in {\Bbb N}$ there are $x_1,\dots,x_n \in C$ such that $\dd( \conv\{x_1,\dots,x_k\},  \conv\{x_{k+1},\dots,x_n\}) \geq \varepsilon$ for all $k=1,\dots, n-1$;
\item[$(\mu_2)$] the supremum of the $\varepsilon>0$ such that there are $\varepsilon$-separated dyadic trees of arbitrary height;
\item[$(\mu_3)$] $=\Delta(C^{\mathcal U})$, for ${\mathcal U}$ a free ultrafilter on ${\Bbb N}$;
\item[$(\mu_4)$] $=\gamma(C^{\mathcal U})$, for ${\mathcal U}$ a free ultrafilter on ${\Bbb N}$;
\item[$(\mu_5)$] the infimum of the $\varepsilon>0$ such that $\Dz(C,\varepsilon) <\omega$;
\item[$(\mu_6)$] the infimum of the $\varepsilon>0$ such that $C$ supports a convex bounded $\varepsilon$-uniformly convex function.
\end{itemize}
Then $\mu_1 \leq \mu_2 \leq 2\mu_3 \leq 2\mu_4 \leq 2\mu_1$ and $\mu_4 \leq 2\mu_5 \leq 2\mu_6 \leq 2\mu_2$.
\end{theo}

\noindent
{\bf Proof.} We will label the steps of the proof by the couple of numbers associated to the inequality.\\
(1-2) If  $\varepsilon < \mu_1$, the separation between convex hulls applied to $2^n$ elements allows the construction of  a 
$\varepsilon$-separated dyadic trees of height $n$. Therefore  $\mu_2 \geq \mu_1$.\\
(2-3) If $\varepsilon < M_2$ then $\Delta(C^{\mathcal U}) \geq \varepsilon/2$. Indeed,  $C^{\mathcal U}$ contains an infinite $\varepsilon$-separated dyadic tree $T$, therefore any nonempty slice of $T$ cannot be covered by finitely many balls of radius less than
$\varepsilon/2$.\\
(3-4) By \cite[Proposition 6.1]{CPR}, $\mbox{Dent}(A) \leq \gamma(A)$, therefore $ \Delta(C^{\mathcal U}) \leq \gamma(C^{\mathcal U})$.\\
(4-1) Let $\varepsilon < \mu_4$. Then there is $x \in \overline{C}^{w^*}$ which is at distance greater than $\varepsilon$ from $X$. Following Oja's proof of James theorem \cite[Theorem 3.132]{banach}, it is posible to find an infinite sequence $(x_n)$ with convex separation greater than $\varepsilon$. Finite representativity gives arbitrarily large sequences in $X$ with the same separation, thus $\mu_1 \geq \mu_4$.\\
(4-5) If $\varepsilon>M_5$ then there is a finite sequence of sets $C=C_1 \supset C_2 \supset \dots \supset C_n$ given by the $\varepsilon$-dentability process. Taking weak$^*$ closures in the bidual, we have
$$ \overline{C}^{w^*} = ( \overline{C_1}^{w^*} \setminus  \overline{C_2}^{w^*}) \cup \dots \cup 
( \overline{C_{n-1}}^{w^*} \setminus  \overline{C_n}^{w^*}) \cup  \overline{C_n}^{w^*}.$$
Now, take any $x \in \overline{C}^{w^*} $ that belongs to one of those sets. The $w^*$-open slice separating $x$ from the smaller set, say $\overline{C_{k+1}}^{w^*}$ ($\emptyset$ for the last set) in the difference is contained in the $w^*$-closure of a slice of $C_k$ not meeting $C_{k+1}$ which has diameter less than $2\varepsilon$ (Lancien's midpoint argument). Since $w^*$-closures does not increase the diameter, we have $d(x,X) \leq 2\varepsilon$. The argument actually implies 
$\gamma(C) \leq 2\varepsilon$. However, we can apply it to the sequence of sets in $X^{\mathcal U}$
$$ C^{\mathcal U} = C_1^{\mathcal U} \supset C_2^{\mathcal U}  \supset \dots \supset C_n^{\mathcal U} $$
which has the same slice-separation property by Lemma \ref{ultra-eps}.\\
(5-6) If $\varepsilon>\mu_6$, there is a bounded convex and $\varepsilon$-uniformly convex function $f$ that, without loss of generality, we may suppose lower semicontinuous. By Proposition \ref{criterio-slice},  any slice of the set 
$\{x \in C: f(x) \leq a\}$
not meeting the set $\{x \in C: f(x) \leq a+\delta\} $ has diameter less than $\varepsilon$. A judicious arranging of these sets shows that $C$ is $\varepsilon$-finitely dentable. Thus $\mu_5 \leq \mu_6$.\\
(6-2) Take $\varepsilon > \mu_2$. Then the $\varepsilon$-separated dyadic trees are uniformly bounded in height. By Theorem \ref{tree}, that implies the existence of $\varepsilon'$-uniformly convex function for every $\varepsilon'>\varepsilon$. Thus $\mu_6 \leq \mu_2$.
\findemo

\begin{rema}
The equivalence between $\mu_3$ and $\mu_4$ is both a local and a quantitative version of the well know statement saying that super-RNP is the same that super-reflexivity. Let us point out that some other relations between the quantities $\mu_i$ for $i=1,\dots,6$ can be established and so improving the equivalence constants. For instance $\mu_2 \leq 2\mu_5$, which is somehow straightforward, or $\mu_6 \leq \mu_5$ as a consequence of Proposition \ref{bush}.
\end{rema}

We will need the following estimation of the distance to the points added by the closure with respect to the topology induced by a norming subspace of the dual.

\begin{lema}\label{lema-cierre}
Let $X$ a Banach space and $F \subset X^*$ an $1$-norming subspace. Then for any bounded convex $A \subset X$ and any $\varepsilon>\gamma(A)$ we have 
$$ \overline{A}^{\sigma(X,F)} \subset A + 2\varepsilon B_X. $$
\end{lema}

\noindent
{\bf Proof.} By \cite[Proposition 3.59]{PVJV}), $ \overline{A}^{w^*} \subset A + 2\varepsilon B_{X^{**}}$. The linear map $p:X^{**} \rightarrow F^*$ defined by
$p(x^{**})=x^{**}|_{F}$ has norm $1$ and satisfies $p( \overline{A}^{w^*} ) = \overline{A}^{\sigma(F^*,F)}$. We may identify 
$p(X)=X$ isometrically into $F^*$ and so we have $ \overline{A}^{\sigma(F^*,F)} \subset A + 2\varepsilon B_{F^*}$. Therefore $ \overline{A}^{\sigma(X,F)} \subset A + 2\varepsilon B_X$ as wished.\findemo\\

We will need the following result that appears as a fact inside the proof of \cite[Theorem 3.1]{KT}. The $1$-norming subspace
$(X^*)^{\mathcal U} \subset (X^{\mathcal U})^* $ will play an important role.

\begin{lema}\label{lematu}
For any $(x^*_n) \in (X^*)^{\mathcal U} $ and $(a_n) \in \conv(A)^{\mathcal U}$, there is 
$(b_n) \in \conv(A^{\mathcal U})$ such that $ \langle (x^*_n), (a_n) \rangle \leq  \langle (x^*_n), (b_n) \rangle$.
\end{lema}

Among the quantities given by Theorem \ref{quantify-swc} only $\mu_4$ does not requiere convexity, so we can propose it as a natural measure of super weak noncompactness. The following is a quantitative version (in terms of $\mu_6$) of  \cite[Theorem 3.1]{KT} establishing that the super weak compactness is stable by closed convex hulls.
Note that the measure of super weak noncompactness introduced by K.Tu in  \cite{KT} is different from ours and so our result is not equivalent to \cite[Theorem 4.2]{KT}.

\begin{theo}
Let $A \subset X$ be a bounded subset and ${\mathcal U}$ a free ultrafilter. Then
$$ \gamma(\conv(A)^{\mathcal U}) \leq 4 \gamma(A^{\mathcal U}). $$
\end{theo}

\noindent
{\bf Proof.} Consider $F=(X^*)^{\mathcal U}$ which is an $1$-norming subspace of $(X^{\mathcal U})^*$. Take 
$\varepsilon>\gamma(A^{\mathcal U})$. By Lemma \ref{lema-cierre},
$$ \overline{\conv}^{\sigma(X,F)}(A^{\mathcal U}) \subset \conv(A^{\mathcal U}) + 2\varepsilon B_{X^{\mathcal U}}.$$
We claim that $(\conv(A))^{\mathcal U} \subset \conv(A^{\mathcal U}) + 2\varepsilon B_{X^{\mathcal U}}$. If it is not the case, then we could separate a point $(\conv(A))^{\mathcal U}$ from $\overline{\conv}^{\sigma(X,F)}(A^{\mathcal U})$ by a functional from $F$. That leads to a contradiction with Lemma \ref{lematu}. Now, we have 
$$ \gamma( (\conv(A))^{\mathcal U} ) \leq 2 \gamma(A^{\mathcal U}) +2\varepsilon $$
which implies the statement.\findemo

\section{A new glance at Enflo's theorem}

Let us show how Enflo's theorem follows from our results.

\begin{theo}[Enflo \cite{Enflo}]
Let $X$ be a super-reflexive Banach space. Then $X$ has an equivalent uniformly convex norm.
\end{theo}

\noindent
{\bf Proof.} The unit ball $B_X$ endowed with the weak topology is SWC. Therefore, for every $\varepsilon>0$,
there is a bounded convex 
$\varepsilon$-uniformly convex function defined on $B_X$ by Theorem \ref{main-swc}.
Now, by Theorem \ref{main-renorming}, there is an equivalent norm $\| \cdot \|_\varepsilon$ on $X$ whose square is an 
$\varepsilon$-uniformly convex function on $B_X$. Without loss of generality, we may assume that $\| \cdot \| \leq \| \cdot \|_\varepsilon \leq 2 \| \cdot \|$. The series $ |\!|\!|  \cdot |\!|\!|^2 =\sum_{n=1}^\infty 2^{-n} \| \cdot \|^2_{1/n} $ defines an equivalent uniformly convex norm.\findemo\\

Enflo's original proof of the uniformly convex renorming of super-reflexive Banach spaces has remain practically unchanged in books, see \cite[pages 438-442]{banach} for instance. We believe that the reason is that the proof is difficult to follow from a geometrical point of view. One of the original aims of this paper was to cast some light on the renorming of super-reflexive spaces. Since the geometrical ideas are now diluted along this paper, we would like to offer to the interested reader a more direct pathway to Enflo's theorem in several steps.

\begin{itemize}

\item From the usual definition of super-reflexivity with finite representation, it is easy to prove that the unit ball $B_X$ of a super-reflexive space has the finite tree property, that is, given $\varepsilon>0$, the uniform boundedness in height of all the $\varepsilon$-separated dyadic trees \cite{James1}.

\item The maximal height of an $\varepsilon$-separated tree with root $x \in B_X$ is an $\varepsilon$-uniformly concave function $h(x)$. This is the main idea in the proof of Theorem \ref{tree}. Note that this function is also symmetric.

\item $g(x)=3^{-h(x)}$ is a symmetric $\varepsilon$-uniformly convex function taking values in $[0,1]$. This comes from Lemma \ref{quasi-3} and is just an arithmetical fact.

\begin{figure}\label{pprince}
\centerline{\includegraphics[width=5cm]{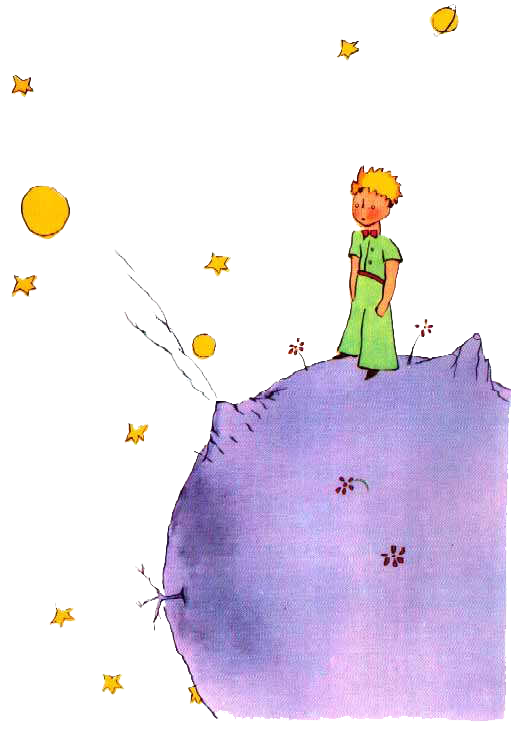}}
\caption{Antoine de Saint-Exupéry's ``Le Petit Prince''}
\end{figure}

\item $f=\breve{g}$ is convex, symmetric and $3\varepsilon$-uniformly convex. The key idea is that $f(x)$ is computed with the values of $g(y)$ with $\|y-x\|<\varepsilon$. The technical details can be carried out as in the proof of Theorem \ref{main-convexify}, which relies on Proposition \ref{bounded-reduction}. Nevertheless, the idea is very intuitive: Planet Earth is a non-convex $\varepsilon$-uniformly convex radial body for $\varepsilon=800 \,\mbox{km}$ at most (see Figure 3). That implies you do not need the Rocky Mountains neither the Himalayas to compute the convex hull over France.

\item Let $f_n$ be the function given in the previous steep for $\varepsilon=1/n$. The function 
$$F(x)=\|x\|+\sum_{n=1}^\infty 2^{-n} f_n(x) $$ 
is uniformly convex, symmetric, Lipschitz on the balls $rB_X$ for $0<r<1$ and it attains a strong minimum at $0$. Moreover,  elementary computations can show that $F(0) \leq 1/17$ and $F(x) \geq 1$ for $x \in S_X$.

\item The set $B=\{x: F(x) \leq \inf F+1/2 \}$ is the unit ball of an equivalent uniformly convex norm $ |\!|\!| \cdot  |\!|\!|$. The idea is to use the Lipschitz property of $F$ to show that for any $\delta>0$, there is $\lambda(\delta) \in (0,1)$ such that 
$$ \{x: F(x) \leq \inf F+1/2 -\delta \} \subset \lambda(\delta) B. $$
Therefore, if $ |\!|\!| x  |\!|\!| =  |\!|\!| y  |\!|\!| =1$ and $\| x-y \| \geq \varepsilon$ then  $F(x)=F(y)=\inf F+1/2$ and 
$F(\frac{x+y}{2}) \leq \inf F +1/2-\delta$ for some $\delta=\delta(\varepsilon)>0$, and thus $  |\!|\!|  \frac{x+y}{2}  |\!|\!|  \leq \lambda(\delta)$.

\end{itemize}

\section*{Acknowledgements}

 The authors are deeply grateful to the referee for the thorough reading of the manuscript and the valuable suggestions made.
 We want also to thank P. Vandaële for his technical help.


\begin{thebibliography}{99}


\bibitem{AP}{\sc D. Azé, J. P. Penot},
Uniformly convex and uniformly smooth convex functions,
{\it Annales de la faculté des sciences de Toulouse}, 6e série, tome 4, no 4 (1995), 705--730.

\bibitem{BL}{\sc Y. Benyamini, J. Lindenstrauss}, 
{\it ~Geometric Nonlinear Functional Analysis.
Vol.~1}, American Mathematical Society Colloquium Publications, 48. American Mathematical Society, Providence, RI, 2000.

\bibitem{BGHV}
{\sc J. M. Borwein, A. Guirao, P. H\'{a}jek, V. Vanderwerff},
{Uniformly convex functions on Banach spaces},
{\it Proc. Amer. Math. Soc.} 137 (2009), 1081--1091.
\newblock \href {https://doi.org/10.1090/S0002-9939-08-09630-5}
  {\path{doi: 10.1090/S0002-9939-08-09630-5}}.

 

\bibitem{BV}
{\sc J. M. Borwein, V. Vanderwerff},
{\it Convex functions: constructions, characterizations and counterexamples},
Encyclopedia of Mathematics and its Applications, 109. Cambridge University Press, Cambridge, 2010.


\bibitem{BV2}
{\sc J. M. Borwein, V. Vanderwerff},
Constructions of Uniformly Convex Functions,
{\it Canad. Math. Bull.} Vol. 55 (4), 2012 pp. 697--707.
\newblock \href {https://dx.doi.org/10.4153/CMB-2011-049-2}
  {\path{doi: 10.4153/CMB-2011-049-2}}.

\bibitem{bourgin}{\sc R. D. Bourgin}, 
{\it Geometric aspects of convex sets with the {R}adon-{N}ikod\'ym property}, 
Lecture Notes in Mathematics, vol. 993.
Springer-Verlag, Berlin (1983).

\bibitem{CMR}{\sc B. Cascales, W. Marciszewski, M. Raja},
{Distance to spaces of continuous functions}. 
{\it Topology Appl.} 153 (2006), no. 13, 2303--2319. 
\newblock \href {https://doi.org/10.1016/j.topol.2005.07.002}
  {\path{doi:10.1016/j.topol.2005.07.002}}.


\bibitem{CPR}{\sc
B. Cascales, A. Pérez, M. Raja},
Radon--Nikodým indexes and measures of weak noncompactness. 
{\it J. Funct. Anal.} 267 (2014), no. 10, 3830--3858. 
\newblock \href {https://dx.doi.org/10.1016/j.jfa.2014.09.015}
  {\path{doi: 10.1016/j.jfa.2014.09.015}}.

\bibitem{cepe}{\sc  M. Cepedello-Boiso}
{Approximation of Lipschitz functions by $\Delta$-convex functions in Banach spaces},
{\it Israel J. Math.} 106 (1998) 269--284.

\bibitem{Cheng}{\sc L. Cheng, Q. Cheng, B. Wang, W. Zhang},
{On super-weakly compact sets and uniformly convexifiable sets},  
{\it Studia Math.} 199 (2010), no. 2, 145--169. 
\newblock \href {https://dx.doi.org/10.4064/sm199-2-2}
  {\path{doi: 10.4064/sm199-2-2}}.


\bibitem{Cla}{\sc J. A. Clarkson},
{`Uniformly convex spaces'},
{\it Trans. Amer. Math. Soc.}  40  (1936),  no. 3, 396--414.

\bibitem{DGZ}{\sc R. Deville, G. Godefroy, V. Zizler}
{\it Smoothness and renormings in Banach spaces},
Pitman Monographs and Surveys in Pure and Applied Mathematics, 64. Longman Scientific \& Technical,
Harlow, 1993.


\bibitem{Enflo}{\sc P. Enflo},
{Banach spaces which can be given an equivalent uniformly convex norm},
{\it Israel J. Math.} 13 (1972), 281--288.


\bibitem{banach} {\sc M. Fabian, P. Habala, P. H\'ajek, V. Montesinos, V. Zizler},
{\it Banach Space Theory. The Basis for Linear and Nonlinear Analysis},
CMS Books in Mathematics,
Springer, New York, 2011.

\bibitem{FHMZ}{\sc M. Fabian, P. Hájek, V. Montesinos, V. Zizler}, 
A quantitative version of Krein's theorem, 
{\it Rev. Mat. Iberoamericana 21} (2005) 237--248.

\bibitem{GL}{\sc L. C. García-Lirola, M. Raja},
Maps with the Radon--Nikodym property. 
{\it Set-Valued Var. Anal.} 26 (2018), 77--93. 
\newblock \href {https://dx.doi.org/10.1007/s11228-017-0428-5}
  {\path{doi: 10.1007/s11228-017-0428-5}}.

\bibitem{GoKi}{\sc K. Goebel, Kazimierz, W. A. Kirk},
{\it Topics in metric fixed point theory.}
Cambridge Studies in Advanced Mathematics, 28. 
Cambridge University Press, Cambridge, 1990.

\bibitem{granero}{\sc A. S. Granero}, 
An extension of the Krein-Šmulian theorem, 
{\it Rev. Mat. Iberoam.} 22(1) (2006), 93--110.

\bibitem{PVJV}{\sc P. H\'ajek, V. Montesinos, J. Vanderwerff, V. Zizler},
{\it Biorthogonal Systems in Banach Spaces},
CMS Books in Mathematics,
Springer, New York, 2018.

\bibitem{James1}{\sc R. C. James},
{Super-reflexive Banach spaces},
{\it Canad. J. Math.} {24} (1972), 896 -- 904.

\bibitem{JL}{\sc W. B. Johnson, J. Lindenstrauss},
{Basic concepts in the geometry of Banach spaces},
{\it Handbook of the Geometry of Banach spaces}
Vol. 1, W.B. Johnson and J. Lindenstrauss editors,
Elsevier, Amsterdam (2001), 1-- 84.

\bibitem{LR}{\sc G. Lancien, M. Raja},
{Nonlinear aspects of super weakly compact sets},
{\it Preprint} (2020)
\newblock \href {https://arxiv.org/abs/2003.01030}
  {\path{arXiv: 2003.01030}}.


\bibitem{LePo}{\sc E. S. Levitin, B. T. Polyak},
{Convergence of minimizing of sequences in the conditional-extremum problem},
{\it Dokl. Akad. Nauk SSSR} 168 (1966), 997--1000.


\bibitem{LT}
{\sc J.~Lindenstrauss, L.~Tzafriri}, 
{\em Classical {B}anach spaces. {II}},
  vol.~97 of Ergebnisse der Mathematik und ihrer Grenzgebiete [Results in
  Mathematics and Related Areas], Springer-Verlag, Berlin-New York, 1979.


\bibitem{Pisier}{\sc G. Pisier},
{Martingales with values in uniformly convex spaces},
{\it Israel J. Math.} 20 (1975), no. 3-4, 326--350.

\bibitem{raja}{\sc M. Raja},
{Finitely dentable functions, operators and sets},
{\it J. Convex Anal.} {15} (2008), 219--233.

\bibitem{raja2}{\sc M. Raja},
{Super WCG Banach spaces},
{\it J. Math. Anal. Appl.} 439 (2016), no. 1, 183--196.
\newblock \href {https://dx.doi.org/10.1016/j.jmaa.2016.02.057}
  {\path{doi: 10.1016/j.jmaa.2016.02.057}}.


\bibitem{vlad1}{\sc A. A. Vladimirov, Yu. E. Nesterov, Yu. N. Chekanov}, 
{On uniformly convex functionals}, {\it Vestnik Moskov. Univ. Ser. XV Vy\v{c}isl Mat. Kibernet.} No. 3 (1978), 12--23 [Russian].

\bibitem{vlad2}{\sc A. A. Vladimirov, Yu. E. Nesterov, Yu. N. Chekanov}, 
{On uniformly quasi- convex functional}, {\it Vestnik Moskov. Univ. Ser. XV Vy\v{c}isl Mat. Kibernet.} No. 4 (1978),
18-27 


\bibitem{KT}{\sc K. Tu}, {Convexification of super weakly compact sets and measure of super weak noncompactness}, 
to appear {\it Proc. Amer. Math. Soc.} (2021).
\newblock \href {http://dx.doi.org/10.1090/proc/15393}
  {\path{doi: 10.1090/proc/15393}}.


\bibitem{Za1}{\sc C. Z\v{a}linescu},
{On uniformly convex functions}, 
{\it J. Math. Anal. Appl.} 95 (1983), no. 2, 344--374.
\newblock \href {https://dx.doi.org/10.1016/0022-247X(83)90112-9}
  {\path{doi: 10.1016/0022-247X(83)90112-9}}.


\bibitem{Za2}{\sc C. Z\v{a}linescu},
{\it Convex analysis in general vector spaces}, 
World Scientific Publishing Co. Inc., River Edge, NJ, 2002. 

\end{thebibliography}
\end{document}